\documentclass{scrartcl}

\usepackage[top=2cm, bottom=2.25cm, left=2.5cm, right=2.5cm]{geometry}
\usepackage[T1]{fontenc}
\usepackage[utf8]{inputenc}
\usepackage{amsmath, amsthm, amssymb}
\usepackage{dsfont}
\usepackage{graphicx}
\usepackage{color}
\usepackage[format=plain]{caption}
\usepackage{capt-of}
\usepackage{subcaption}
\usepackage{booktabs}
\usepackage{multicol}
\usepackage{bm}
\usepackage{dlfltxbcodetips} 
\usepackage{tabularx}
\usepackage{csquotes}
\usepackage{marvosym}

\newcommand{\norm}[1]{\left\| #1 \right\|}  
\renewcommand{\d}{\,\mathrm{d}} 
\newcommand{\N}{\mathbb{N}}  

\newcommand{\R}{\mathbb{R}}

\newcommand{\dist}{\operatorname{dist}}

\newcommand{\diag}{\operatorname{diag}}

\newcommand{\eps}{\varepsilon}
\renewcommand{\phi}{\varphi}
\newcommand{\ul}{\underline}
\newcommand{\ol}{\overline}
\newcommand{\tr}{\operatorname{tr}}

\newcommand{\Cov}{\operatorname{Cov}}
\newcommand{\Var}{\operatorname{Var}}
\newcommand{\dimx}{d_\mathrm{x}}
\newcommand{\dimy}{d_\mathrm{y}}
\newcommand{\dimp}{d_\mathrm{\theta}}
\newcommand{\Rpsd}{\R_{\mathrm{psd}}}
\newcommand{\Rpd}{\R_{\mathrm{pd}}}
\newcommand{\X}{\mathcal{X}}
\newcommand{\Y}{\mathcal{Y}}
\newcommand{\xtot}{x^{\mathrm{tot}}}
\newcommand{\ytot}{y^{\mathrm{tot}}}
\newcommand{\txtot}{\tilde{x}^{\mathrm{tot}}}
\newcommand{\tytot}{\tilde{y}^{\mathrm{tot}}}
\newcommand{\ntot}{n^{\mathrm{tot}}}
\newcommand{\nsam}{n^{\mathrm{sam}}}
\newcommand{\lin}{\mathrm{lin}}
\newcommand{\sam}{\mathrm{sam}}
\newcommand{\lse}{\widehat{\theta}}
\newcommand{\flin}{f_{\ol{\theta}}^{\lin}}
\newcommand{\init}{\texttt{init}}
\newcommand{\fed}{\texttt{fed}}
\newcommand{\oed}{\texttt{oed}}
\newcommand{\tot}{\texttt{tot}}

\numberwithin{equation}{section}

\newtheorem{thm}{Theorem}

\newtheorem{ass}[thm]{Assumption}
\theoremstyle{definition} \newtheorem{algo}{Algorithm}
\theoremstyle{definition} 
\theoremstyle{definition}

\title{Sequential optimal experimental design for vapor-liquid equilibrium modeling}

\author{Martin Bubel$^1$, Jochen Schmid$^1$, Volodymyr Kozachynskyi$^2$,\\ 
Erik Esche$^2$, Michael Bortz$^1$\\  
\small $^1$Fraunhofer Institute for Industrial Mathematics (ITWM), 67663 Kaiserslautern, Germany\\ 
\small $^2$Technische Universit\"at Berlin, Process Dynamics and Operations Group, 10623 Berlin, Germany\\
\small jochen.schmid@itwm.fraunhofer.de}  

\date{}

\begin{document}

\maketitle

\begin{abstract}
\small{ \noindent 
We propose a general methodology of sequential locally optimal design of experiments for explicit or implicit  nonlinear models, as they abound in chemical engineering and, in particular, in vapor-liquid equilibrium modeling. As a sequential design method, our method iteratively alternates between performing experiments, updating parameter estimates, and computing new experiments. Specifically, our sequential design method computes a whole batch of new experiments in each iteration and this batch of new experiments is designed in a two-stage locally optimal manner. In essence, this means that in every iteration the combined information content of the newly proposed experiments and of the already performed experiments is maximized. 
In order to solve these two-stage locally optimal design problems, a recent and efficient adaptive discretization algorithm is used.
We demonstrate the benefits of the proposed methodology on the example of of the parameter estimation for the non-random two-liquid model for narrow azeotropic vapor-liquid equilibria. As it turns out, our sequential optimal design method requires substantially fewer experiments than traditional factorial design to achieve the same model precision and prediction quality. Consequently, our method can contribute to a substantially reduced experimental effort in vapor-liquid equilibrium modeling and beyond. 
}
\end{abstract}


\section{Introduction}

In order to reliably simulate and optimize engineering processes, one needs sufficiently precise 
mathematical models describing the functional relationships between relevant input and output quantities of the considered process. In many cases, these models come in the form of parametric models, that is, their predictions $f(x,\theta)$ for the output quantity $y$ of interest not only depend on the input quantity values $x$ but also on additional model parameters $\theta$. A priori, 
the true value of these model parameters is unknown and therefore has to be estimated from appropriate experimental data. Specifically, at suitably selected input points $x_1, \dots, x_n$ one measures  the values $y_1, \dots, y_n$ of the output quantity and then adjusts the model parameters $\theta$ in such a way that the discrepancy between the model's predictions $f(x_i,\theta)$ and the measurements $y_i$ becomes minimal in the least-squares sense.  
\smallskip

Intuitively, it is clear that in general the quality of the thus obtained parameter estimate strongly depends on the selected experimental points $x_1, \dots, x_n$ or, in other words, on where one measures and on how often one measures. It is also clear that in practice 
one can often perform only a fairly limited number $n$ of experiments, just because experiments in a lab are typically time- and cost-intensive. Consequently, it is of crucial importance to carefully plan and select the actually performed experiments. In particular, it is important that (i) the parameter estimate and the resulting predictions based on the selected experiments $x_1, \dots, x_n$ are sufficiently precise and that (ii) the number $n$ of experiments is affordably small. Coming up with experimental plans that meet these criteria is the core task of experimental design. 
\smallskip


In the extensive literature on experimental design, two main approaches can be distinguished, namely factorial experimental design on the one hand and optimal experimental design on the other hand. In factorial design, one spreads the experiments in a uniform and grid-like manner in the entire input space and, as a consequence, ends up with a relatively large number of experiments already for moderate input space dimensions. In optimal design, by constrast, one tries to 
make experiments only at places in the input space where enough information can be gained about the model parameters to be estimated. 
Specifically, one tries to distribute the experimental points such that the information content -- and hence the precision -- of the corresponding parameter estimate becomes optimal. 
As a consequence, with optimal design one can expect to obtain the same parameter and prediction quality as with factorial design, but with substantially fewer experiments. 
And this effect 
can be expected to be all the more pronounced, the more nonlinear the considered model $f$ is, as a function of the inputs $x$ or of the model parameters $\theta$. In light of this, optimal experimental design is a promising approach especially for vapor-liquid equilibrium modeling with its intrinsically 
nonlinear models for activity coefficients, for instance. 
\smallskip 


In the present paper, we propose a general methodology of sequential locally optimal experimental design~\cite{Ch59, Si, FeLe, PrPa} for general nonlinear models, 
and these nonlinear models can be defined either by explicit or by implicit algebraic equations. As a sequential design method, our method iteratively alternates between carrying out experiments, estimating the parameters, and designing new experiments. Specifically, the sequential design method proposed here is characterized by the following key features, which also distinguish it from the sequential design method recently proposed in~\cite{DuAt21}. 
\smallskip

A first essential feature of our method is that, in every iteration, it computes a whole batch of new experiments, the maximal batch size being specified by the user. In contrast, the design methodology of~\cite{DuAt21} computes new experiments one by one (except for the initial iteration). 
We decided for batches of new experiments because this will generally lead to fewer alternations between the experimental and the computational phase of the procedure and thus to a smoother workflow in practice. In particular, this is beneficial when each experimental phase requires extensive set-up or ramp-up times or when the experiments and the computations are performed by different people at different institutions. 
\smallskip

A second important feature of our method is that, in every iteration, the batch of new experiments is computed in a locally optimal manner, namely by solving a two-stage locally optimal experimental design problem. 
Specifically, the new experiments are designed such that the combined information content of the already performed experiments (stage one) and of the new experiments to be performed next (stage two) becomes maximal. As usual in two-stage locally optimal design, the combined information content is measured in terms of a suitable design criterion $\Psi$ applied to the information matrix of a suitably linearized version of the model $f$ at the combined experimental design. Also, the linearization of $f$, in every iteration, is performed about the current best parameter estimate, namely, 
about the least-squares parameter estimate based on all experimental data collected until the current iteration. As the design criterion $\Psi$, in turn, one can choose each of the usual standard criteria provided only it is differentiable. In particular, this includes the A- and the D-criterion, while it excludes essentially only the E-criterion. 
\smallskip

A third distinctive feature of our method is that it solves the occurring two-stage locally optimal experimental design problems by means of the recent adaptive discretization algorithm from~\cite{YaBiTa13}, which is specifically tailored to such design problems. And, in fact, this algorithm exhibits very good convergence properties, both from a theoretical and from a practical point of view. 
We terminate our iterative alternation between experiments, parameter estimation, and locally optimal experimental design, if the total number of already performed and of newly proposed experiments exceeds a user-specified upper bound or if the newly proposed experiments do not differ substantially anymore from already performed ones.
\smallskip

In order to demonstrate the benefits of the proposed sequential design methodology for practical use, we apply it to the estimation of non-random two-liquid parameters. As is well-known, these parameters describe the binary interactions in the liquid phase of a multi-component mixture and therefore they play an essential role in vapor-liquid equilibrium modeling. 
It is also well-known that the estimation of these parameters 
is generally a challenging task, especially for azeotropic systems. In fact, the corresponding sum of squared errors often features flat local minima, multiple local minima, or even multiple global minima~\cite{Ta79, JaAs83, We23}.  
\smallskip

In the present paper, we specifically consider the binary system consisting of propanol and propyl acetate and we measure the vapor composition and the temperature of this system at vapor-liquid equilibrium for different liquid compositions and pressures. It is known that this system exhibits narrow azeotropic behavior and therefore the estimation of its non-random two-liquid parameters is already a rather challenging task. Applying our sequential design method to this system, we end up with a locally optimal experimental design consisting of $15$ experiments. And, as it turns out, these $15$  optimally designed experiments lead to the same prediction quality as $27$ experiments from a traditional factorial experimental design. 
So, compared to traditional factorial design, our sequential locally optimal design method substantially reduces the experimental effort -- and thus helps the experimenter to save substantial amounts of time and money -- without impairing the model prediction quality. 
Additionally, our method also helps the experimenter in deciding what to do in case there still is some experimental budget left, after all the experiments proposed by a first run of our sequential design method have been performed. Indeed, by its very definition, our method always indicates where 
making a batch of further experiments is most promising and where it is not.  
\smallskip

Apart from the papers~\cite{DuAt20, DuAt21}, we are not aware of any other works that apply optimal design techniques 
to the estimation of thermodynamic property parameters. 
It should be noticed that in the mentioned  papers~\cite{DuAt20, DuAt21} 
no real lab experiments are performed, but simulated data are used instead, based on literature values for the non-random two-liquid parameters to be estimated. In sharp contrast to that, 
we perform actual lab experiments and only rely on the corresponding real experimental data. In particular, we do not rely on literature data or on literature values for the considered parameters at all.  


\section{Setting and preliminaries} \label{sec:setting-and-preliminaries}

In this section, we specify the setting and the basic terminology used throughout the paper. We also collect those facts from parameter estimation and experimental design that are relevant for the sequential locally optimal experimental design method proposed in this paper. 

\subsection{Setting and basic terminology} \label{sec:setting-and-terminology}

In the entire paper, we assume that the considered functional relationship between input quantities $x$ and output quantities $y$ can be described by the predictions $f(x,\theta)$ of a parametric model
\begin{align}
f: \X \times \Theta \to \Y,
\end{align}
where $\X \subset \R^{\dimx}$ and $\Theta \subset \R^{\dimp}$ are the sets of admissible input values and of admissible model parameter values respectively, while $\Y := \R^{\dimy}$ is a superset of the set of possibly observable output values. 
\smallskip

Such a parametric model $f$ can be linear or nonlinear and, moreover, it can be explicit or implicit. 
As usual, $f$ is called \emph{linear} iff it is affine-linear w.r.t.~$\theta$, that is, 
\begin{align} \label{eq:linear-model}
f(x,\theta) = c(x) + J(x)\theta \qquad ((x,\theta) \in \X \times \R^{\dimp})
\end{align}
for some arbitrary functions $c: \X \to \R^{\dimy}$ and $J: \X \to \R^{\dimy \times \dimp}$. 
In particular, a linear model $f$ can be arbitrarily nonlinear w.r.t.~the inputs $x \in \X$. Simple examples of linear models are furnished by the polynomial functions in $x$ with the monomial coefficients making up the model parameter vector $\theta$.
Additionally, $f$ is called \emph{explicit} iff it is explicitly defined by a closed-form expression for $f(x,\theta)$ for every $(x,\theta) \in \X \times \Theta$. And correspondingly, $f$ is called \emph{implicit} iff one has no closed-form expression for $f(x,\theta)$ but instead $f(x,\theta)$ for every $(x,\theta) \in \X \times \Theta$ is implicitly defined as the solution $y = h(s(x,\theta),x,\theta)$ of a system of algebraic equations of the form
\begin{align} \label{eq:implicit-model-definition}
g(s,x,\theta) = 0 \qquad \text{and} \qquad y = h(s,x,\theta)
\end{align}
with explicitly defined functions $g$ and $h$. In this system of equations, the first equation~(\ref{eq:implicit-model-definition}.a) determines the state $s = s(x,\theta)$ as a function of $(x,\theta)$, while the second equation~(\ref{eq:implicit-model-definition}.b) determines the output quantity $y = h(s(x,\theta),x,\theta)$ as a function of the state. 
\smallskip

Simple examples of explicit nonlinear models from chemical engineering are given by the Antoine model for the saturation vapor pressure of a pure substance or by the non-random two-liquid model for the activity coefficients of a mixture of substances. Similarly, a simple chemical engineering example of an implicit nonlinear model arises in the prediction 
of the bubble-point temperature of a mixture of substances. With this particular example of an implicit nonlinear model, we will  illustrate the general methodology proposed in this paper. See Section~\ref{sec:application} below.
\smallskip

We point out, though, that the proposed general methodology is by no means restricted to thermodynamic models but rather can be applied in a wide range of circumstances, to explicit and implicit nonlinear models alike. All that needs to be satisfied for the methodology to apply are a few rather mild regularity and representability assumptions on the considered model $f$. Specifically, the model $f$ needs to be sufficiently regular (Assumption~\ref{ass:regularity}) and sufficiently representative (Assumption~\ref{ass:observation-representability}) in the sense that, up to additive random measurement errors, all possible output measurement values are explained 
by the model's predictions with one specific model parameter value $\theta^*$. 
Spelled out in detail, these assumptions take the following form.

\begin{ass} \label{ass:regularity}
\begin{itemize}
\item[(i)] $\X$ is a compact subset of $\R^{\dimx}$ and $\ul{\X}$ is a compact subset of $\X$
\item[(ii)] $\Theta$ is a compact subset of $\R^{\dimp}$ that is equal to the closure of its interior, in short, $\ol{\Theta^\circ} = \Theta$
\item[(iii)] $\Theta \ni \theta \mapsto f(x,\theta)$ is continuously differentiable for every $x \in \X$ and $\X \ni x \mapsto D_\theta f(x,\ol{\theta})$ is continuous for every $\ol{\theta} \in \Theta$. As usual, $D_{\theta}f(x,\overline{\theta})$ denotes the derivative (Jacobi matrix) of $f(x,\cdot)$  at $\overline{\theta}$.
\end{itemize}
\end{ass}

Conventionally, the subset $\ul{\X}$ of the input space $\X$ from the assumption above is called the \emph{experimental design space}. It consists of all potential experiments, 
that is, all input points $x \in \X$ at which the experimenter can measure the output quantity. 
We point out that the experimental design space $\ul{\X}$ can be discrete (finite) or continuous (infinite). In particular, $\ul{\X}$ can be equal to $\X$, but in most experimental design algorithms $\ul{\X}$ is chosen to be finite for the sake of compuational feasibility~\cite{FeLe, PrPa, YaBiTa13, DuAt21, DuAt21}. We also point out that the assumed regularity conditions are very mild and usually easy to verify. Indeed, for an explicit model one has only to investigate the explicit model parametrization $f$, and for an implicit model one has only to investigate the implicit model-defining functions $g$ and $h$ and invoke the implicit function theorem. 

\begin{ass} \label{ass:observation-representability}
A parameter value $\theta^* \in \Theta$ exists such that for every $n \in \N$ and every set 
of input values $x_1, \dots, x_n \in \ul{\X}$, the corresponding measured values $y_1, \dots, y_n \in \Y$ of the output quantity are given by the predictions $f(x_1,\theta^*), \dots, f(x_n,\theta^*)$ of the model $f(\cdot,\theta^*)$ up to independent normally distributed measurement errors $\epsilon_1, \dots, \epsilon_n$, that is, 
\begin{align} \label{eq:representability}
y_i = f(x_i,\theta^*) + \epsilon_i \qquad (i \in \{1,\dots,n\})
\end{align}
where $\epsilon_1, \dots, \epsilon_n$ are realizations of independent and normally distributed measurement errors $\eps_1, \dots, \eps_n$ having mean $0$ and a known positive definite covariance matrix $\varsigma \in \R^{\dimy \times \dimy}_{\mathrm{pd}}$.
\end{ass}

Any value $\theta^*$ as in the above assumption is called (the) \emph{true value of the model parameter} and, accordingly, the corresponding model $f(\cdot,\theta^*)$ is called (the) \emph{true model}. It is important to notice that the assumption does not require such a true parameter value to be known to the modeler. In fact, quite to the contrary, it is always unknown in practice and therefore has to be estimated. All the assumption requires is that such a true (but unknown) parameter value exists and that the additive random measurement errors are independent and normally distributed. 
See~\cite{BaWa} (Section~1.3.1) for a critical discussion and justification of these standard assumptions in parameter estimation~\cite{BaWa, SeWi} and in experimental design~\cite{FeLe, Mo, PrPa}.
\smallskip

With the basic assumptions set out, we now formally introduce the notion of (unweighted and weighted) experimental designs 
along with the pertinent notation that will be continually used throughout the paper. 
An \emph{experimental design} on $\ul{\X}$ simply is a tuple 
\begin{align} \label{eq:def-unweighted-design}
\tilde{x} = (x_1,\dots,x_n) \in \ul{\X}^n
\end{align}
consisting of a natural number $n$ of input points $x_1, \dots, x_n \in \ul{\X}$ which are called the \emph{experimental points of $\tilde{x}$}. 
Also, $n$ is called the \emph{size of the experimental design $\tilde{x}$}. 
It should be noticed that the experimental points of a given experimental design $\tilde{x}$ are not required to be pairwise distinct. 
If they are, then $\tilde{x}$ is called an \emph{experimental design without replications (repetitions)}, and otherwise $\tilde{x}$ is called an \emph{experimental design with replications (repetitions)}. 
Since performing actual experiments in a lab is typically expensive, it is often important to take into account experiments $\tilde{x}^- = (x^{-}_1, \dots, x^{-}_{n^-})$ that have already be performed at an earlier stage and combine them with newly performed experiments $\tilde{x}^+ = (x^{+}_1, \dots, x^{+}_{n^+})$. We will denote the resulting \emph{combined (or refined) experimental design} by
\begin{align} \label{eq:combined-design}
\tilde{x}^- \,\&\, \tilde{x}^+ := (x^{-}_1, \dots, x^{-}_{n^-}, x^{+}_1, \dots, x^{+}_{n^+}) \in \ul{\X}^{n^- + n^+}.
\end{align}
In order to distinguish experimental designs in the sense of~\eqref{eq:def-unweighted-design} from weighted experimental designs defined next, we will sometimes refer to the former as \emph{unweighted experimental designs}. 
\smallskip

If $\tilde{x} = (x_1,\dots,x_n)$ is a given unweighted experimental design, it is often useful to aggregate the repeated experimental points, to indicate the number of replications, and then to normalize by the size $n$ of the design. In this manner, one can associate to every unweighted experimental design $\tilde{x} = (x_1,\dots,x_n)$ the corresponding weighted experimental design
\begin{align} \label{eq:def-weighted-design-corresponding-to-unweighted-design}
\xi_{\tilde{x}} := \big((w_x,x))_{x\in\ul{\X}}
\qquad \text{with} \qquad w_x := r_x/n,
\end{align}
where $r_x \in \{0,\dots, n\}$ for every potential experiment $x \in \ul{\X}$ indicates how often it occurs in the design $\tilde{x}$. Sometimes, it is further useful to allow for general real weights $w_x \in [0,1]$ instead of the specific rational weights in~\eqref{eq:def-weighted-design-corresponding-to-unweighted-design}. Specifically, for a finite design space $\ul{\X}$, it is sometimes useful to consider tuples $((w_x,x))_{x\in\ul{\X}}$ of general weight-experiment pairs. Such tuples are called \emph{weighted experimental designs} on $\ul{\X}$ and are usually denoted by the letter $\xi$ in the literature. Accordingly, the set of all weighted experimental designs on a finite design space $\ul{\X}$ is denoted by
\begin{align}
\Xi(\ul{\X}) 
:= \bigg\{\big((w_x,x))_{x\in\ul{\X}}: w_x \in [0,1] \text{ for every } x \in \ul{\X} 
\text{ and } \sum_{x \in \ul{\X}} w_x = 1\bigg\}. 
\end{align}

In parameter estimation and in experimental design alike, information matrices play a fundamental role. As usual, the \emph{one-point information matrix of $f$ at a design point $x \in \ul{\X}$} for the reference parameter $\ol{\theta}$ is defined as
\begin{align} \label{eq:one-point-information-matrix}
m_f(x,\ol{\theta}) := D_{\theta} f(x,\ol{\theta})^\top \varsigma^{-1} D_{\theta} f(x,\ol{\theta})
\in \Rpsd^{\dimp\times\dimp},
\end{align} 
where $\Rpsd^{\dimp\times\dimp}$ denotes the set of all positive semidefinite $\dimp\times\dimp$ matrices. Analogously, the \emph{information matrix of an experimental design} -- either an unweighted experimental design $\tilde{x} = (x_1, \dots, x_n)$ or a weighted experimental design $\xi = ((w_x,x))_{x\in\ul{\X}} \in \Xi(\ul{\X})$ -- is  defined as the sum or, respectively, as the weighted sum of the corresponding one-point information matrices. In formulas, 
\begin{gather} \label{eq:information-matrix-for-designs}
M_f(\tilde{x},\ol{\theta}) := \sum_{i=1}^n D_{\theta} f(x_i,\ol{\theta})^\top \varsigma^{-1} D_{\theta} f(x_i,\ol{\theta})
\in \Rpsd^{\dimp\times\dimp},\\
M_f(\xi,\ol{\theta}) := \sum_{x\in\ul{\X}} w_x D_{\theta} f(x,\ol{\theta})^\top \varsigma^{-1} D_{\theta} f(x,\ol{\theta})\in \Rpsd^{\dimp\times\dimp}.
\end{gather} 

An extension of the tilde notation~\eqref{eq:def-unweighted-design} for experimental designs to vector- or matrix-valued functions on $\X$ or $\X \times \Theta$ is very convenient both in parameter estimation and in experimental design. Specifically, if $a: \X \to \mathcal{V}$ and $b: \X \times \Theta \to \mathcal{V}$ are given functions from $\X$ or $\X \times \Theta$ to some space $\mathcal{V}$ of real column vectors or of real matrices and if $\tilde{x} = (x_1, \dots, x_n)$ is a given experimental design of size $n$, then we will repeatedly use the shorthand notations 
\begin{align}
\tilde{a}(\tilde{x}) 
:=
\begin{pmatrix}
a(x_1)\\
\vdots \\
a(x_n)
\end{pmatrix} \in \mathcal{V}^n
\qquad \text{and} \qquad
\tilde{b}(\tilde{x}, \theta) 
:=
\begin{pmatrix}
b(x_1, \theta)\\
\vdots \\
b(x_n, \theta)
\end{pmatrix} \in \mathcal{V}^n
\end{align}
to denote the outcome of applying $a$ or $b$ simultaneously to all experimental points of $\tilde{x}$. With this notation, Assumption~\ref{ass:observation-representability} can be concisely reformulated as follows: for every experimental design $\tilde{x} = (x_1,\dots,x_n)$ on $\ul{\X}$, the corresponding experimental results $\tilde{y} = (y_1^\top, \dots, y_n^\top)^\top$ are realizations of the $n\dimy$-dimensional random variable
\begin{align} \label{eq:random-observation-model}
\tilde{y}(\tilde{x}) := \tilde{f}(\tilde{x},\theta^*) + \tilde{\eps}(\tilde{x}),
\end{align}
where $\tilde{\eps}(\tilde{x}) = (\eps_1^\top, \dots, \eps_n^\top)^\top$ is an $n\dimy$-dimensional normally distributed random variable with expectation value $0$ and block-diagonal covariance matrix
\begin{align}
\tilde{\varsigma} 
:=
\diag(\varsigma, \dots, \varsigma) 
= 
\begin{pmatrix}
\varsigma & \dotsb & 0\\
\vdots & \ddots & \vdots\\
0 & \dotsb & \varsigma
\end{pmatrix} \in \Rpd^{n\dimy \times n\dimy}.
\end{align} 
Similarly, applying the above tilde notation to the function $J: \X \times \Theta \to \R^{\dimy \times \dimp}$ defined by $J(x,\ol{\theta}) := D_\theta f(x,\ol{\theta})$, the information matrix of an experimental design $\tilde{x} = (x_1, \dots, x_n)$ can be rewritten in the form 
\begin{align} \label{eq:information-matrix-alternative-representation}
M_f(\tilde{x}, \ol{\theta}) = \tilde{J}(\tilde{x},\ol{\theta})^\top \tilde{\varsigma}^{-1} \tilde{J}(\tilde{x},\ol{\theta}).
\end{align}
It is clear from this representation that the information matrix $M_f(\tilde{x},\ol{\theta})$ is positive definite if and only if the matrix $\tilde{J}(\tilde{x},\ol{\theta}) \in \R^{n\dimy\times\dimp}$ has a trivial nullspace. In particular, for the information matrix $M_f(\tilde{x},\ol{\theta})$ of $\tilde{x} = (x_1,\dots,x_n)$ to be invertible, it is necessary that 
\begin{align} \label{eq:necessary-condition-for-invertibility}
l\dimy \ge \dimp
\qquad (l:=\#\{x_1,\dots,x_n\}) 
\end{align}
or, in other words, that the number $l$ of different experimental points in $\tilde{x}$ is at least as large as $\dimp/\dimy$, the number of model parameters divided by the number of outputs. As is well-known, the invertibility of $M_f(\tilde{x},\ol{\theta})$ plays an essential role in the estimability of the model parameters from the experimental design $\tilde{x}$, see~\eqref{eq:exact-lse-for-linear-models} and~\eqref{eq:approximate-lse-for-nonlinear-models} below and~\cite{PrPa} (Section 7.2). 

\subsection{Central facts from parameter estimation} \label{sec:parameter-estimation}

Since the true value $\theta^*$ of the model parameter (Assumption~\ref{ass:observation-representability}) is unknown, it needs to be estimated. 
A standard way of doing so is least-squares estimation, that is, one measures the values $y_1, \dots, y_n \in \Y$ of the output quantity at certain experimental points $x_1, \dots, x_n \in \ul{\X}$ and then chooses the model parameters $\theta$ such that the model's predictions $f(x_1,\theta), \dots, f(x_n,\theta)$ deviate as little as possible from the observed measurements $y_1, \dots, y_n$. Specifically, one chooses the model parameter $\theta \in \Theta$ as a solution to the least-squares estimation problem
\begin{align} \label{eq:lse-problem}
\min_{\theta \in \Theta} \sum_{i=1}^n \big(f(x_i,\theta)-y_i\big)^\top \varsigma^{-1} \big(f(x_i,\theta)-y_i\big),
\end{align} 
where $\varsigma \in \Rpd^{\dimp\times\dimp}$ is the covariance matrix of the measurement errors (Assumption~\ref{ass:observation-representability}). Any solution to~\eqref{eq:lse-problem} is called a \emph{least-squares estimate for $f$ based on the experimental design $\tilde{x} = (x_1,\dots,x_n)$ with corresponding  experimental results $\tilde{y} = (y_1^\top, \dots, y_n^\top)^\top$}. Since $\Theta$ is assumed to be compact and $f(x_i,\cdot)$ is continuous (Assumption~\ref{ass:regularity}), the existence of a least-squares estimate is guaranteed for arbitrary experimental data $\tilde{x}$ and $\tilde{y}$. If the model is further assumed to be globally least-squares estimable for a given experimental design $\tilde{x}$, then uniqueness of the least-squares estimate for $\tilde{x}$ 
can also be established (Section~7.2 of~\cite{PrPa}). We do not need such global estimability assumptions, however,  and will simply write $\lse_f(\tilde{x},\tilde{y})$ for any of the potentially multiple least-squares estimates for $f$ based on $\tilde{x}$ and $\tilde{y}$.  

\subsubsection{Computing least-squares estimates}

In the special case of linear models, one has a closed-form expression for the least-squares estimates. Specifically, for $f$ as in~\eqref{eq:linear-model}, the least-squares estimate is unique and given by the formula
\begin{align} \label{eq:exact-lse-for-linear-models}
\lse_f(\tilde{x},\tilde{y}) = \big( \tilde{J}(\tilde{x})^\top \tilde{\varsigma}^{-1} \tilde{J}(\tilde{x}) \big)^{-1} \tilde{J}(\tilde{x})^\top \tilde{\varsigma}^{-1} \big( \tilde{y} - \tilde{c}(\tilde{x}) \big),
\end{align}
provided that the occurring inverse of $\tilde{J}(\tilde{x})^\top \tilde{\varsigma}^{-1} \tilde{J}(\tilde{x})$ exists, see~\eqref{eq:necessary-condition-for-invertibility}. If it does not, then one can replace it in~\eqref{eq:exact-lse-for-linear-models} by the Moore-Penrose pseudoinverse to get the unique minimal-norm least-squares estimate of $f$ for $\tilde{x}$ and $\tilde{y}$. 
As is well-known, the pseudoinverse can be computed by means of the singular value decomposition or by various regularizations. See~\cite{NoWr} (Section~10.2) and~\cite{GoVa96} (Chapter 5) for an overview of concrete methods to compute least-squares estimates for linear models.
\smallskip

In the general case of nonlinear models, by contrast, one usually has no closed-form expression for the least-squares estimates anymore and therefore one has to resort to suitable approximations. A standard approximation is based on  linearizing $f$ around a suitable reference parameter value $\ol{\theta} \in \Theta$. Specifically, in the least-squares problem~\eqref{eq:lse-problem} one replaces the general nonlinear model $f$ by the linearized model $\flin$ defined by
\begin{align} \label{eq:linearized-model}
f_{\ol{\theta}}^{\lin}(x,\theta) := f(x,\ol{\theta}) + J(x,\ol{\theta})(\theta-\ol{\theta})
\qquad ((x,\theta) \in \X \times \R^{\dimp})
\end{align}
with $J(x,\ol{\theta}) := D_{\theta}f(x,\ol{\theta})$. Applying then~\eqref{eq:exact-lse-for-linear-models}, one obtains the approximate formula
\begin{align} \label{eq:approximate-lse-for-nonlinear-models}
\lse_f(\tilde{x},\tilde{y}) \approx \lse_{\flin}(\tilde{x},\tilde{y}) 
= \ol{\theta} + \big( \tilde{J}(\tilde{x},\ol{\theta})^\top \tilde{\varsigma}^{-1} \tilde{J}(\tilde{x},\ol{\theta}) \big)^{-1} \tilde{J}(\tilde{x},\ol{\theta})^\top \tilde{\varsigma}^{-1} \big( \tilde{y} - \tilde{f}(\tilde{x},\ol{\theta}) \big)
\end{align}
provided that the occurring inverse of $\tilde{J}(\tilde{x},\ol{\theta})^\top \tilde{\varsigma}^{-1} \tilde{J}(\tilde{x},\ol{\theta}) = M_f(\tilde{x},\ol{\theta})$ exists, see~\eqref{eq:information-matrix-alternative-representation} and~\eqref{eq:necessary-condition-for-invertibility} above. If one iterates this approximation~\eqref{eq:approximate-lse-for-nonlinear-models} starting from an initial $\ol{\theta}^0$ and if in each iteration one chooses an appropriate step size, then one arrives at the Gau\ss-Newton method. As is well-known, this iterative method converges to a critical point of the 
sum-of-squares objective function in~\eqref{eq:lse-problem}, provided that the matrices $\tilde{J}(\tilde{x},\theta)^\top \tilde{\varsigma}^{-1} \tilde{J}(\tilde{x},\theta)$ are invertible for all $\theta$ in the region of interest 
(Theorem~10.1 of~\cite{NoWr}). If this uniform invertibility is not satisfied, one can use the Levenberg-Marquardt method instead (Theorem~10.3 of~\cite{NoWr}). We refer to~\cite{NoWr} (Section~10.3) and to~\cite{SeWi} (Chapter~14) for a an extensive overview of methods to compute least-squares estimates for nonlinear models.

\subsubsection{Covariance of least-squares estimates}

Since the experimental results $\tilde{y}$ for any given experimental design $\tilde{x}$ are subject to random measurement errors, the same is true for the corresponding least-squares parameter estimates $\lse_f(\tilde{x},\tilde{y})$. In other words, the least-squares estimates for any given experimental design $\tilde{x}$ are uncertain. A standard measure to quantify this uncertainty is the covariance matrix
\begin{align}
\Cov \lse_f(\tilde{x},\tilde{y}(\tilde{x}))
\in \Rpsd^{\dimp \times \dimp}
\end{align} 
of the least-squares estimate for $\tilde{x}$ and the random observations $\tilde{y}(\tilde{x})$ from~\eqref{eq:random-observation-model} (Assumption~\ref{ass:observation-representability}). 
In the case of linear models, by~\eqref{eq:exact-lse-for-linear-models} and~\eqref{eq:information-matrix-alternative-representation} one has the exact identity
\begin{align} \label{eq:exact-covariance}
\Cov \lse_f(\tilde{x},\tilde{y}(\tilde{x})) = M_f(\tilde{x},\ol{\theta})^{-1}
\end{align}
for arbitrarily chosen reference parameter values $\ol{\theta} \in \Theta$. Indeed, the information matrix $M_f(\tilde{x},\ol{\theta})$ actually does not at all depend on $\ol{\theta}$ for linear models~\eqref{eq:linear-model}. 
In the case of nonlinear models, by~\eqref{eq:approximate-lse-for-nonlinear-models} and~\eqref{eq:information-matrix-alternative-representation} one has at least the approximate identity
\begin{align} \label{eq:approximate-covariance}
\Cov \lse_f(\tilde{x},\tilde{y}(\tilde{x})) \approx M_f(\tilde{x},\ol{\theta})^{-1}
\end{align}
for suitably chosen reference parameter values $\ol{\theta} \in \Theta$. As opposed to linear models, the information matrix $M_f(\tilde{x},\ol{\theta})$ of nonlinear models generally strongly depends on the reference parameter value $\ol{\theta}$. In general, the approximation~\eqref{eq:approximate-covariance} will be the better, the closer $\ol{\theta}$ is to true parameter value $\theta^*$. A rigorous asymptotic justification of this 
is provided by the classical asymptotic normality results for least-squares estimation (Section~12.2.3 of~\cite{SeWi} or Section~3.1.3 of~\cite{PrPa}). 

\subsection{Central facts from experimental design} \label{sec:experimental-design}

Since by~\eqref{eq:exact-covariance} and~\eqref{eq:approximate-covariance} the uncertainty of the least-squares estimate $\lse_f(\tilde{x},\tilde{y}(\tilde{x}))$ depends on the underlying experimental design $\tilde{x}$, the prediction uncertainty of the estimated model $f(\cdot, \lse_f(\tilde{x},\tilde{y}(\tilde{x})))$ depends on the experimental design as well. And so, in order to get a sufficiently precise estimated model that is sufficiently certain about its predictions, one has to carefully choose the experimental design $\tilde{x}$ on which the estimation is based. In experimental design, one tries to come up with experimental designs $\tilde{x} = (x_1,\dots,x_n)$ that on the one hand have a fesible size $n$ and on the other hand lead to an estimated model that is sufficiently certain about its predictions. Commonly, this uncertainty is quantified in terms of a suitable scalar function
\begin{align}
\Psi(\Cov \lse_f(\tilde{x},\tilde{y}(\tilde{x}))^{-1})
\end{align}
of the inverse of the estimate's covariance. Such functions $\Psi$ are called \emph{design criteria} and for the proposed methodology to work, we have to impose only the following assumptions on $\Psi$ which are standard in  convex design theory~\cite{Ki74, Si, Pu, AtDoTo, FeLe, PrPa}. A reader not familiar with that theory can simply always think of the particularly important special case where $\Psi$ is the D- or the A-criterion, see ~\eqref{eq:def-D-and-A-criterion} below. 

\begin{ass} \label{ass:design-criterion}
$\Psi: \Rpsd^{\dimp \times \dimp} \to \R \cup \{\infty\}$ is convex, antitonic, lower semicontinuos function whose restriction to $\Rpd^{\dimp \times \dimp}$ takes only finite values and is continuously differentiable. 
\end{ass}

All standard design criteria $\Psi_p$ except for the E-criterion $\Psi_{\infty}$ satisfy the above assumptions and thus the methodology proposed here can be applied with all standard design criteria except for the E-criterion. In particular, the assumptions are satisfied for the D-criterion $\Psi_0$ and for the A-criterion $\Psi_1$, which are defined as follows:
\begin{align} \label{eq:def-D-and-A-criterion}
\Psi_0(M) := \ln \det(M^{-1}) 
\qquad \text{and} \qquad
\Psi_1(M) := \tr(M^{-1})
\end{align}
for $M \in \Rpd^{\dimp \times \dimp}$ and $\Psi_p(M) := \infty$ for $M \in \Rpsd^{\dimp \times \dimp}\setminus \Rpd^{\dimp \times \dimp}$. 
%
%
See~\cite{ScSeBo24} for precise definitions of the notions used in Assumption~\ref{ass:design-criterion} and of the standard design criteria $\Psi_p$ as well as for detailed verifications. 
\smallskip

In the extensive literature on experimental design, two main approaches can be distinguished, namely factorial design~\cite{Fi, BoHu, Mo} on the one hand and optimal design~\cite{KiWo59, Fe, Ki74, Si, Pu, AtDoTo, FeLe, PrPa} on the other hand. In essence, these approaches differ in their primary goal and in their generality and scope of applicability. 

\subsubsection{Factorial experimental design} \label{sec:FED}

In factorial experimental design, the primary goal is to find a factorial experimental design $\tilde{x}^*$, that is, a geometrically regular design whose experimental points constitute either a full grid in $\ul{\X}$ or a regularly formed fraction of a full grid in $\ul{\X}$. 
Accordingly, one speaks of full factorial designs in the first case and of fractional factorial designs in the second case. 
In contrast to optimal experimental design, the optimality of~$\tilde{x}^*$ plays only 
a secondary role here. 
As it turns out, for certain fairly special models $f$ and for cuboidal design spaces 
\begin{align} \label{eq:cuboidal-design-space}
\ul{\X} = [a_1,b_1] \times \dotsb \times [a_{\dimx},b_{\dimx}],
\end{align}
suitably chosen factorial designs are also optimal (w.r.t.~suitable design criteria $\Psi$). Specifically, for linear models $f$ that are also linear or at least multilinear in the inputs $x$, the standard full factorial designs consisting of the corner points of~\eqref{eq:cuboidal-design-space} are A-, D-, and E-optimal (Section~3.3.4 of~\cite{Ba}). 
Similarly, for linear models $f$ that are quadratic in the inputs $x$, fractional factorial designs can be found~\cite{Ko62} that are D-optimal. It should be noted, however, that these designs are weighted designs and their support is fairly large even for small input dimensions (Section~2.2.2 of~\cite{Ba}). It should also be noted that, unlike the specialized designs from~\cite{Ko62}, the standard factorial designs for linear models that are quadratic in the inputs are not D-optimal (\cite{NaGoMi70} or Section~3.4.7 of~\cite{Ba}).
As soon as the model $f$ becomes more general (higher-order in $x$ or nonlinear in $\theta$), 
factorial designs will no longer be optimal, however. 
In fact, this can already be seen from third-order polynomial regression in one input dimension (Section~2.2.1 of~\cite{Ba} 
or Section~2.6 of~\cite{FeLe}) or in two input dimensions (Section~2.2.2 of~\cite{Ba}). 
So, for applications described by higher-order linear models or by nonlinear models, factorial design is not the method of choice and optimal experimental design is recommended instead (Section~11.4.4 of~\cite{Mo}). Similarly, optimal experimental design is recommended for models on non-cuboidal input spaces $\X$ (Section~11.4.4 of~\cite{Mo}). 

\subsubsection{Optimal experimental design} \label{sec:OED}

In optimal experimental design, the primary goal is to find an optimal experimental design $\tilde{x}^*$ w.r.t.~the chosen design criterion $\Psi$, that is, a solution to the optimization problem
\begin{align} \label{eq:OED}
\min_{\tilde{x} \in \ul{\X}^n} \Psi(\Cov \lse_f(\tilde{x},\tilde{y}(\tilde{x}))^{-1}),
\end{align} 
where the design size $n$ can be fixed or allowed to vary in some given range. 
In contrast to factorial design, the geometric regularity of $\tilde{x}^*$ plays no role here. And in fact, by the remarks in the previous paragraph, optimal experimental designs are generally not factorial. 
As has been pointed out above, except for 
linear models $f$, one usually has no exact expressions for the estimate's covariance and hence for the objective function in~\eqref{eq:OED}. So, for nonlinear models $f$, one has to work with appropriate approximations of the covariance. In many cases, the linearization-based approximation~\eqref{eq:approximate-covariance} is used, which turns the original optimal experimental design problem~\eqref{eq:OED} into the substitute problem 
\begin{align} \label{eq:local-OED}
\min_{\tilde{x} \in \ul{\X}^n} \Psi(M_f(\tilde{x},\ol{\theta}))
\end{align}
with some appropriate reference parameter value $\ol{\theta} \in \Theta$ and with a fixed or a variable design size $n$. 
It should be noticed that~\eqref{eq:local-OED} is a nonlinear optimization problem which is discrete, continuous, or discrete-continuous depending on whether the design space $\ul{\X}$ is discrete or continuous and on whether $n$ is fixed or not. Conventionally,  problems like~\eqref{eq:local-OED} are called locally optimal experimental design problems~\cite{Ch53} and, accordingly, any approach based on substituting~\eqref{eq:OED} by~\eqref{eq:local-OED} is referred to as locally optimal design. It is important to notice that the linearization-based approximation~\eqref{eq:approximate-covariance} is by no means restricted to locally optimal design. In fact, quite to the contrary, the vast majority of approaches to optimal experimental design for nonlinear models relies, in one way or another, 
on the approximation~\eqref{eq:approximate-covariance}. 
In particular, this is true for average optimal design and minimax optimal design, 
which are defined by the problems
\begin{align}
\min_{\tilde{x} \in \ul{\X}^n} \int_{\Theta} \Psi(M_f(\tilde{x},\theta)) \d\pi(\theta)
\qquad \text{and} \qquad
\min_{\tilde{x} \in \ul{\X}^n} \max_{\theta \in \Theta} \Psi(M_f(\tilde{x},\theta))
\end{align}
with some prior probability distribution $\pi$ on $\Theta$. Additionally, this is true for the class of sequential design methods \cite{Ch59, Si, FeLe, PrPa},
to which the method proposed here belongs. We refer to~\cite{FeLe, PrPa} for comprehensive expositions of these linearization-based approaches to optimal experimental design for nonlinear models and to~\cite{FrMa08} for a concise overview of locally optimal experimental design specifically tailored to dynamic nonlinear models. We also refer to~\cite{HaWa85, ViGa07, MuPa19} for approaches that -- unlike the first-order approaches mentioned so far -- propose higher-order approximations of the covariance or of the related confidence regions of the least-squares estimate.

\section{Sequential locally optimal experimental design} \label{sec:sequential-experimental-design}

In this section, we formally introduce our sequential locally optimal experimental design method. In rough terms, the method repeatedly solves locally optimal experimental design problems of the kind 
\begin{align} \label{eq:local-OED-1}
\min_{\tilde{x} \in \ul{\X}^n} \Psi(M_f(\tilde{x},\ol{\theta}))
\end{align} 
with successively improved reference parameter values $\ol{\theta}$. Specifically, in every iteration, $\ol{\theta}$ is chosen to be the current best parameter estimate, based on all experiments performed until then, and the  locally optimal experimental design problem with that value of $\ol{\theta}$ is solved by means of the recent adaptive discretization algorithm from~\cite{YaBiTa13}. We continue this iteration until the total number of performed experiments reaches some 
upper bound or the newly proposed experiments do not differ substantially anymore from already performed experiments. As is usual~\cite{FeLe, PrPa, YaBiTa13} 
for concrete experimental design algorithms, we choose the design space $\ul{\X}$ to be a finite set 
and work with weighted designs for the sake of computational feasibility. 

\begin{ass} \label{ass:design-space-finite-and-information-matrices-computable}
$\ul{\X}$ is a finite set and for every given reference parameter $\ol{\theta} \in \Theta$ 
all the one-point information matrices $m_f(x,\ol{\theta}) = D_{\theta} f(x,\ol{\theta})^\top \varsigma^{-1} D_{\theta} f(x,\ol{\theta})$ for $x \in \ul{\X}$ can be computed sufficiently fast.
\end{ass}

In the case of an explicit model $f$, the computability assumption above poses no restriction. In the case of an implicit model $f$, in turn, it boils down 
to requiring a sufficiently fast numerical solvability of the equation~(\ref{eq:implicit-model-definition}.a) for $s = s(x,\ol{\theta})$ along with a sufficiently fast numerical invertibility of the Jacobi matrices $D_s g(s(x,\ol{\theta}),x,\ol{\theta})$ for all $x \in \ul{\X}$. 

\subsection{Computing locally optimal experimental designs}

At the heart of our sequential experimental design method is the solution of suitable variants of the locally optimal experimental design problem~\eqref{eq:local-OED-1}. In essence, these variants are obtained by means of two modifications. 
\smallskip

As a first modification to~\eqref{eq:local-OED-1}, we extend the search space from the unweighted experimental designs of a given size $n$ to all weighted experimental designs on $\ul{\X}$. In this manner, we arrive at the one-stage locally optimal experimental design problem 
\begin{align} \label{eq:1-stage-local-OED}
\min_{\xi \in \Xi(\ul{\X})} \Psi\big(M_f(\xi,\ol{\theta})\big).
\end{align}
Since the weighted designs $\xi \in \Xi(\ul{\X})$ are uniquely determined by their weights $w_x$ for $x \in \ul{\X}$, 
the design problem~\eqref{eq:1-stage-local-OED} is actually a weight optimization problem. Since, moreover, the design criterion $\Psi$ is convex (Assumption~\ref{ass:design-criterion}), the weight optimization problem~\eqref{eq:1-stage-local-OED} is a convex nonlinear optimization problem and therefore has a much simpler structure than the original general nonlinear 
design problem~\eqref{eq:local-OED-1}. It is precisely this substantial structural simplification that makes the introduction of weighted experimental designs standard in optimal experimental design. In the literature, 
problems like~\eqref{eq:1-stage-local-OED} are referred to as one-stage experimental design problems because their design criterion does not explicitly take into account experiments that have already been performed at an earlier stage. 
Such preliminary experiments, however, are quite often available in experimental practice. 
\smallskip

As a further modification to~\eqref{eq:local-OED-1}, we therefore explicitly incorporate such previous experiments into the design criterion. Specifically, we consider the two-stage locally optimal experimental design problem 
\begin{align} \label{eq:2-stage-local-OED}
\min_{\xi \in \Xi(\ul{\X})} \Psi\big(\alpha M_f(\xi_{\tilde{x}^-},\ol{\theta}) + (1-\alpha) M_f(\xi,\ol{\theta})\big),
\end{align}
where $\tilde{x}^- = (x_1^-, \dots, x_{n^-}^-)$ are experiments that have already been performed at a previous stage, $\xi_{\tilde{x}^-}$ is the corresponding weighted design~\eqref{eq:def-weighted-design-corresponding-to-unweighted-design}, and $\alpha \in [0,1]$ is an importance factor for those previously performed experiments. In the limiting case $\alpha = 0$, the information contained in the already performed experiments is completely ignored in the design of the new experiments, and thus~\eqref{eq:2-stage-local-OED} reduces to~\eqref{eq:1-stage-local-OED} then. In the limiting case $\alpha = 1$, by contrast, the information contained in the new experiments is completely ignored which makes 
the optimization~\eqref{eq:2-stage-local-OED} pointless, of course. If one chooses $\alpha = 0.5$, then the uncertainty induced by the previous experiments and the uncertainty induced by the prospective experiments are both taken into account in equal measure. 
Just like~\eqref{eq:1-stage-local-OED}, the two-stage experimental design problem~\eqref{eq:2-stage-local-OED} is actually a weight optimization problem and, by the convexity of $\Psi$ (Assumption~\ref{ass:design-criterion}), it is a convex continuous optimization problem. 
\smallskip

In order to solve one- or two-stage locally optimal design problems of the kind~\eqref{eq:1-stage-local-OED} or~\eqref{eq:2-stage-local-OED}, one can resort to a large variety of tailor-made 
iterative algorithms from the literature like  the algorithms from~\cite{Fe, Wy70, At73, SiTiTo78, Bo86, Yu11} or \cite{Ha20}, for instance. See also~\cite{AtDoTo, FeLe, PrPa} for a good overview of the pertinent algorithms known until 2013. In this paper, we use the more recent adaptive discretization algorithm from~\cite{YaBiTa13} because of its very good convergence properties, both from a theoretical and from a practical point of view~\cite{YaBiTa13, ScSeBo24, SeScBo21, VaSe21}. 
Instead of solving~\eqref{eq:2-stage-local-OED} in one shot, 
this algorithm repeatedly solves discretized versions of~\eqref{eq:2-stage-local-OED} with successively refined discretizations. 
And since these discretizations are refined adaptively, the discretized problems to be solved in the course of the algorithm remain fairly small, typically. Specifically, we use the following 
version of the algorithm~\cite{ScSeBo24}.

\begin{algo} \label{algo:YBT}
Input: 
an unweighted design $\tilde{x}^-$ on $\ul{\X}$ with an associated importance factor $\alpha \in [0,1)$, 
a reference parameter value $\ol{\theta} \in \Theta$, 
the one-point information matrices $m_f(x,\ol{\theta})$ for all $x \in \ul{\X}$, 
a small finite subset $\ul{\X}^0$ of $\ul{\X}$ such that $\alpha M_f(\xi_{\tilde{x}^-},\ol{\theta}) + (1-\alpha) M_f(\xi,\ol{\theta})$ is invertible for some design $\xi$ on $\ul{\X}^0$, 
and an optimality tolerance $\eps \in (0,\infty)$. 
Initialize $k = 0$. With these inputs, proceed in the following steps.
\begin{itemize}
\item[1.] Compute a solution $\xi^k \in \Xi(\ul{\X}^k)$ of the discretized version
\begin{align} \label{eq:discretized-design-problem}
\min_{\xi \in \Xi(\ul{\X}^k)} \Psi^{(\alpha)}(M_f(\xi,\ol{\theta}))
\end{align}
of the two-stage locally optimal design problem~\eqref{eq:2-stage-local-OED}, where $\Psi^{(\alpha)}$ is the two-stage design criterion defined by $\Psi^{(\alpha)}(M) := \Psi(\alpha M_f(\xi_{\tilde{x}^-},\ol{\theta}) + (1-\alpha)M)$ for all $M \in \Rpsd^{\dimp \times \dimp}$.

\item[2.] Check how far the computed solution $\xi^k$ of the discretized problem~\eqref{eq:discretized-design-problem} is from being optimal for the original problem~\eqref{eq:2-stage-local-OED}. In order to do so, compute a solution $x^k \in \ul{\X}$ of the strongest optimality violator problem
\begin{align} \label{eq:violator-problem}
\min_{x \in \ul{\X}} \psi^{(\alpha)}(M_f(\xi^k,\ol{\theta}), x),
\end{align} 
where $\psi^{(\alpha)}$ is the sensitivity function of $\Psi^{(\alpha)}$ w.r.t.~$M_f(\Xi(\ul{\X}),\ol{\theta})$ (which is a certain directional derivative of $\Psi^{(\alpha)}$). 
\begin{itemize}
\item In case $\psi^{(\alpha)}(M_f(\xi^k,\ol{\theta}), x^k) < -\eps$, refine the discretization by setting $\ul{\X}^{k+1} := \ul{\X}^k \cup \{x^k\}$ and return to Step~1 with $k$ replaced by $k+1$.
\item In the opposite case, return $\xi^k$ and terminate.
\end{itemize}
\end{itemize}
\end{algo}

As is shown in~\cite{ScSeBo24}, the algorithm above terminates after finitely many iterations with an approximately optimal design $\xi^+ \in \Xi(\ul{\X})$ in the sense that
\begin{align} \label{eq:eps-optimal-design-definition}
\Psi\big(\alpha M_f(\xi_{\tilde{x}^-},\ol{\theta}) + (1-\alpha) M_f(\xi^+,\ol{\theta})\big) \le \min_{\xi \in \Xi(\ul{\X})} \Psi\big(\alpha M_f(\xi_{\tilde{x}^-},\ol{\theta}) + (1-\alpha) M_f(\xi,\ol{\theta})\big) + \eps.
\end{align}
In words, this relation says that the design criterion value for $\xi^+$ is guaranteed to exceed the optimal design criterion value by at most $\eps$. It is important to notice that this approximation error bound $\eps$ can be freely set by the user. In the following, we will refer to a design $\xi^+ \in \Xi(\ul{\X})$ that satisfies \eqref{eq:eps-optimal-design-definition} as a \emph{locally $\eps$-optimal weighted design for $f$ on $\ul{\X}$ around $\ol{\theta}$ with $\tilde{x}^-$ as the previous-stage design of importance $\alpha$}. Algorithm~\ref{algo:YBT} is an adaptive discretization algorithm because 
the discretization $\ul{\X}^{k+1}$ is obtained from the previous discretization $\ul{\X}^k$ by an adaptive refinement, namely by adding a strongest violator $x^k$ of the optimality condition for the solution $\xi^k$ of the previous discretized design problem~\eqref{eq:discretized-design-problem}. 
It is because of this adaptive refinement that Algorithm~\ref{algo:YBT} typically does not need many iterations until termination. In particular, the discretized design problems~\eqref{eq:discretized-design-problem} to be solved in the course of the algorithm are typically fairly small and hence easy to solve. After all, these discretized design problems~\eqref{eq:discretized-design-problem} are convex weight optimization problems with the number of weight variables being equal to the size of the discretization $\ul{\X}^k$. Also, by the finiteness of $\ul{\X}$ (Assumption~\ref{ass:design-space-finite-and-information-matrices-computable}), the strongest violator problems~\eqref{eq:violator-problem} can simply be solved by means of enumeration, using the explicit formulas for the sensitivity function $\psi^{(\alpha)}$ of two-stage design criteria $\Psi^{(\alpha)}$ provided in~\cite{ScSeBo24}, for instance. 
\smallskip

As is indicated already by their name, the locally $\eps$-optimal weighted designs computed by Algorithm~\ref{algo:YBT} are weighted experimental designs $\xi^+ = ((w_x^+,x))_{x\in\ul{\X}} \in \Xi(\ul{\X})$. And in most cases,  the weights $w_x^+$ will not be of the form $w_x^+ = r_x^+/n^+$ with a practically realizable number $n^+$ of new experiments and replication numbers 
$r_x^+ \in \{0,\dots,n^+\}$ indicating how often the experimenter should measure at the points $x \in \ul{\X}$. In other words, the weighted experimental designs $\xi^+$ computed by Algorithm~\ref{algo:YBT} usually do not directly correspond 
to an unweighted experimental design $\tilde{x}^+ = (x_1^+, \dots, x_{n^+}^+)$ of a practically realizable size $n^+$. Instead, they first have to be converted to such a practically realizable unweighted design. Standard methods of doing so are based on suitable rounding procedures for the (pseudo)quotas $n^+ w_x^+$ \cite{FeLe} (Section~3.3), $(n^+-l^+)w_x^+$ \cite{Fe} (Section~3.1), or $(n^+-l^+/2)w_x^+$ \cite{PuRi92} (Section~6), where $l^+$ is the number of experimental points in $\xi^+$. In the case where the affordable number $n^+$ of new experiments is large, these rounding procedures produce unweighted designs $\tilde{x}^+$ that are efficient in the sense that the design criterion value of $\xi_{\tilde{x}^+}$ exceeds that of $\xi^+$ only by a little bit. Indeed, for large experimental budgets $n^+$, the efficiency bounds from~\cite{Fe, FeLe, PuRi92} become tight. 
In the practically relevant case of small experimental budgets $n^+$, however, the aforementioned efficiency bounds 
do not say much about the actually achievable efficiency. Additionally, the mentioned efficiency bounds are valid 
only for homogeneous design criteria~\cite{FeLe}, 
but two-stage design criteria -- as the criteria $\Psi^{(\alpha)}$ used here -- are not homogeneous. We therefore propose a different conversion procedure in this paper, namely one that is not based on rounding but instead directly aims at producing efficient unweighted designs $\tilde{x}^+$ from $\xi^+$. Specifically, among the experimental points of $\xi^+$ with non-negligible weights, we simply select $n^+$ points such that the design criterion value becomes as small as possible. 
\smallskip

Summarizing, we arrive at the following two-stage locally optimal experimental design algorithm. Starting from an unweighted design $\tilde{x}^-$ containing experiments from a previous stage, it computes an unweighted experimental design $\tilde{x}^+$ containing the experiments for the next stage. 

\begin{algo} \label{algo:2-stage-locally-optimal-design}
Input: 
an unweighted experimental design $\tilde{x}^-$ on $\ul{\X}$ 
with an associated importance $\alpha \in [0,1)$, 
the one-point information matrices $m_f(x,\ol{\theta})$ for all $x \in \ul{\X}$, 
a reference parameter value $\ol{\theta} \in \Theta$, 
an optimality tolerance $\eps \in (0,\infty)$,
a minimal total weight $\ul{w}^+ \in (0,1]$, 
and a maximal number $\ol{n}^+$ of new experiments.
With these inputs, proceed in the following steps.
\begin{itemize}

\item[1.] Compute a locally $\eps$-optimal weighted design 
\begin{align}
\xi^+ = ((w_x^+,x))_{x\in\ul{\X}}
\end{align}
for $f$ on $\ul{\X}$ around $\ol{\theta}$ with $\tilde{x}^-$ as the previous-stage design of importance $\alpha$. In order to do so, use Algorithm~\ref{algo:YBT}. 

\item[2.] Convert the computed weighted design $\xi^+$ to an unweighted design $\tilde{x}^+$ consisting of at most $\ol{n}^+$ 
experiments. In order to do so, successively sieve out the least important experiment from $\xi^+$ until the total weight of the remaining experiments only just remains above the prescribed threshold $\ul{w}^+$. Write 
\begin{align}
\tilde{\ul{x}}^+ = (\ul{x}_1^+, \dots, \ul{x}_{\ul{n}^+}^+)
\end{align}
for the remaining unweighted design and then compare the number $\ul{n}^+$ of experiments in $\tilde{\ul{x}}^+$ with the prescribed upper bound $\ol{n}^+$. 
\begin{itemize}
\item In case $\ul{n}^+ > \ol{n}^+$, choose an appropriate subdesign $\tilde{x}^+$ of the intermediate design $\tilde{\ul{x}}^+$, consisting of exactly $\ol{n}^+$ experiments. Specifically, among all subdesigns $\tilde{z}^+$ of $\tilde{\ul{x}}^+$ consisting of exactly $\ol{n}^+$ experiments, choose one with minimal design criterion value 
\begin{align}
\Psi\big(\alpha M_f(\xi^-,\ol{\theta}) + (1-\alpha) M_f(\xi_{\tilde{z}^+},\ol{\theta})\big)
\end{align}
\item In the opposite case, simply choose $\tilde{x}^+$ to be equal to the intermediate design $\tilde{\ul{x}}^+$.
\end{itemize}
\end{itemize}
\end{algo}

Clearly, the output $\tilde{x}^+$ of the above algorithm is an unweighted experimental design without replications consisting of at most $\ol{n}^+$ experiments in $\ul{\X}$. In the following, we will refer to it as a \emph{locally $\eps$-optimal design for $f$ on $\ul{\X}$ around $\ol{\theta}$ consisting of at most $\ol{n}^+$ experiments with total weight at least $\ul{w}^+$ and with $\tilde{x}^-$ as the previous-stage design of importance $\alpha$}. 

\subsection{Iterating experiments, parameter estimation, and locally optimal experimental design}

\begin{figure}[!ht]
\centering
\includegraphics[width=0.6\columnwidth]{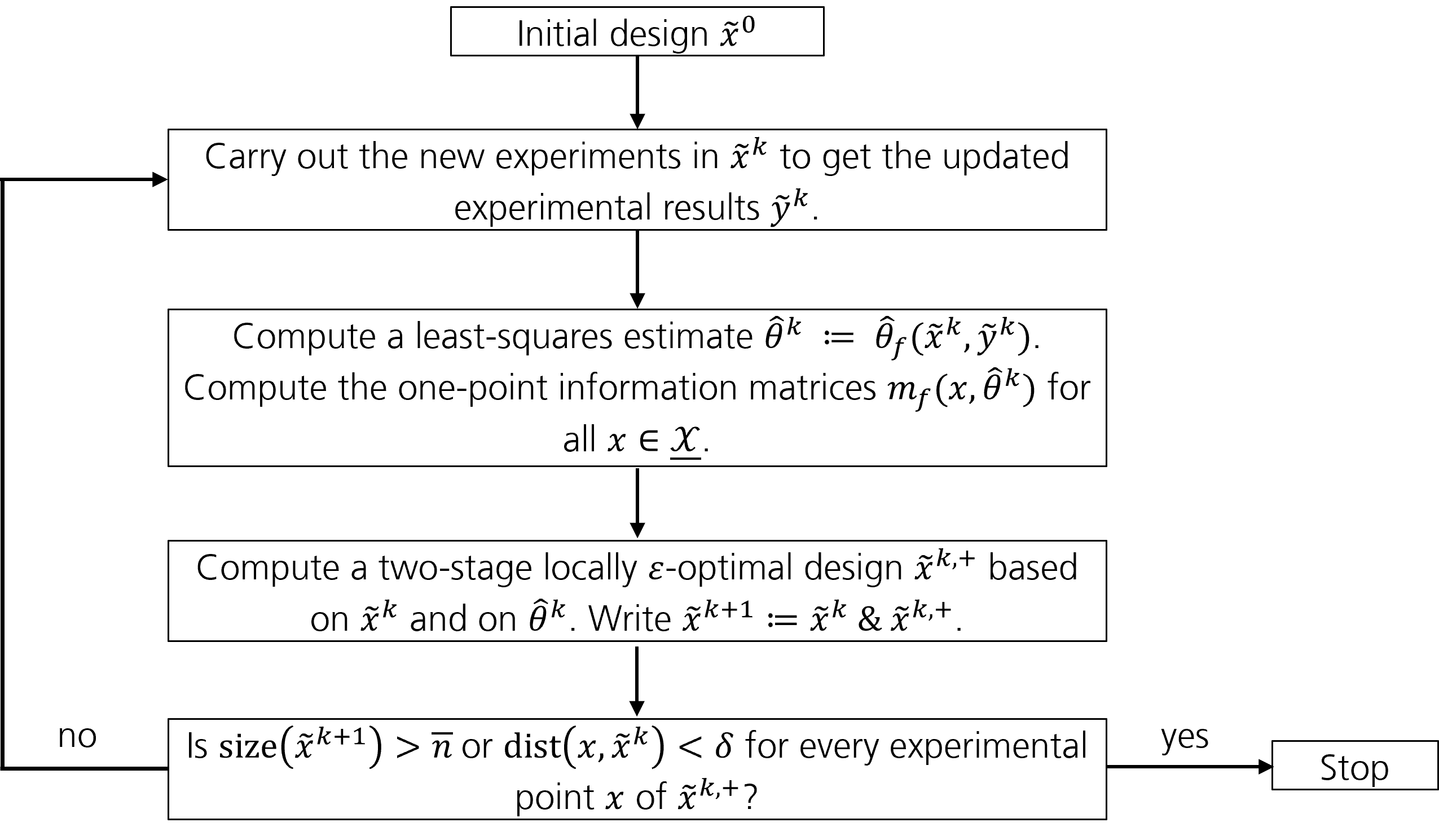}
\caption{Schematic of the proposed sequential locally optimal experimental design method.}
\label{fig:sequential-design-method}
\end{figure}

As has been pointed out above, the quality of the approximation~\eqref{eq:approximate-covariance} generally strongly depends on the specific choice of the reference parameter $\ol{\theta}$ about which the model is linearized. So, as  the locally optimal experimental design problems~\eqref{eq:1-stage-local-OED} and~\eqref{eq:2-stage-local-OED} are built precisely on this approximation, it is essential to come up with good choices for the reference parameter $\ol{\theta}$. 
A very natural and 
established strategy here is sequential design~\cite{Ch59, Si, FeLe, PrPa} where $\ol{\theta}$ is simply always chosen to be the current best parameter estimate $\lse_f(\tilde{x}^-, \tilde{y}^-)$, based on all experimental data $\tilde{x}^-$ and $\tilde{y}^-$ that are available at the current stage. With this estimate $\ol{\theta}$, we then apply our locally optimal experimental design algorithm (Algorithm~\ref{algo:2-stage-locally-optimal-design}) to get the most informative batch $\tilde{x}^+$ of new experiments. We repeat this three-step procedure of performing experiments, estimating the parameters, and computing locally optimal experimental designs until the total number of experiments exceeds a user-specified upper bound $\ol{n}$ or each of the newly proposed experiments differs from an already performed experiment by less than a user-specified progress tolerance $\delta$. We measure this distance between two  experimental points $x = (x_j)_{j\in\{1,\dots,\dimx\}}$ and $x' = (x'_j)_{j\in\{1,\dots,\dimx\}}$ in terms of a suitably scaled maximum norm defined by
\begin{align}
\norm{x-x'}_{\infty} := \max_{j\in\{1,\dots,\dimx\}} |x_j-x'_j|/\lambda_j,
\end{align}
where $\lambda_1, \dots, \lambda_{\dimx}$ are the side lengths of the smallest cuboid (axis-aligned box) containing the design space $\ul{\X}$. 
(In other words, $\lambda_1, \dots, \lambda_{\dimx}$  describe the extension of the design space in the individual input directions.) 
In detail, our sequential design procedure is described in Algorithm~\ref{algo:sequential-optimal-design}. A schematic illustration of the procedure is given in Figure~\ref{fig:sequential-design-method}. 

\begin{algo} \label{algo:sequential-optimal-design}
Input: an initial experimental design $\tilde{x}^0$ on $\ul{\X}$, 
an importance factor $\alpha \in [0,1)$ for existing experiments, 
an optimality tolerance $\eps \in (0,\infty)$,
a minimal total weight $\ul{w}^+ \in (0,1]$, 
a maximal number $\ol{n}^+$ of new experiments per iteration,
a maximal total number $\ol{n}$ of experiments,
and a progress tolerance $\delta \in (0,\infty)$.
Initialize $k = 0$. With these inputs, proceed in the following steps.
\begin{itemize}
\item[1.] Carry out the experiments that have been newly proposed in the previous iteration, that is, measure the output quantity at those experimental points of $\tilde{x}^k = (x_1^k, \dots, x_{n_k}^k)$ that have been newly proposed. Write 
\begin{align}
\tilde{y}^k := ((y_1^k)^\top, \dots, (y_{n_k}^k)^\top)^\top.
\end{align}
for all experimental results collected so far (the newly and the previously collected ones).  

\item[2.] Compute a least-squares estimate $\lse^k	:= \lse_f(\tilde{x}^k, \tilde{y}^k)$ for $f$ based on all experimental data $(\tilde{x}^k, \tilde{y}^k)$ collected so far, and then compute the one-point information matrices $m_f(x,\lse^k)$ for all candidate design points $x \in \ul{\X}$.

\item[3.] Compute a locally $\eps$-optimal design $\tilde{x}^{k,+}$ for $f$ on $\ul{\X}$ around $\lse^k$ consisting of at most $\ol{n}^+$ experiments with total weight at least $\ul{w}^+$ and with the current design $\tilde{x}^k$ as the previous-stage design of importance $\alpha$. In order to do so, use Algorithm~\ref{algo:2-stage-locally-optimal-design}.  Write
\begin{align}
\tilde{x}^{k+1} := \tilde{x}^k \,\&\, \tilde{x}^{k,+}
\end{align}
for the combined design consisting of the current experiments and the newly computed experiments, $n_{k+1} := n_k + n_k^+$ for the size of $\tilde{x}^{k+1}$, and $\dist(x,\tilde{x}^k) := \min_{i\in\{1,\dots,n_k\}} \norm{x-x_i}_{\infty}$ for the distance of any given $x$ from the experimental points of $\tilde{x}^k$. 
\begin{itemize}
\item In case $n_{k+1} > \ol{n}$ or $\dist(x,\tilde{x}^k) < \delta$ for every experimental point $x$ of $\tilde{x}^{k,+}$,  terminate.
\item In the opposite case, return to Step 1 with $k$ replaced by $k+1$.
\end{itemize}
\end{itemize}
\end{algo}

Clearly, if the progress tolerance $\delta$ is chosen small enough, then the progress-based termination criterion of the above algorithm reduces to the condition that each of the newly computed experimental points of $\tilde{x}^{k,+}$ is actually equal to one of the experimental points of the previous design $\tilde{x}^k$. 
Specifically, this is the case if $\delta$ is chosen to be smaller than the mesh size (granularity) 
\begin{align}
\min \big\{ \norm{x-x'}_\infty: x,x' \in \ul{\X} \text{ with } x \ne x' \big\}
\end{align}
of the chosen finite design space $\ul{\X}$ (Assumption~\ref{ass:design-space-finite-and-information-matrices-computable}). 

\subsection{Assessing the quality of experimental designs} \label{sec:quality-measures}

In this section, we introduce the criteria that will be used to assess the quality of a given experimental design. We use three different criteria, and all of them express the quality of an experimental design $\tilde{x} = (x_1, \dots, x_n)$ in terms of the prediction quality of the corresponding estimated model $f(\cdot,\lse_f(\tilde{x},\tilde{y}(\tilde{x})))$, that is, the model that has been trained based on $\tilde{x}$.
\smallskip

As a first quality measure of a given experimental design $\tilde{x}$ (with associated experimental results $\tilde{y}$), we take the sum of squared prediction errors that the trained model $f(\cdot,\lse_f(\tilde{x},\tilde{y}))$ 
makes on all available experimental data $(\txtot,\tytot)$. Specifically, we take the root mean squared errors 
\begin{align} \label{eq:rmse-definition}
\rho_j(\tilde{x}) := \rho_j(\tilde{x},\tilde{y}) := \bigg( \frac{1}{\ntot} \sum_{i=1}^{\ntot} \big( f_j(\xtot_i, \lse_f(\tilde{x},\tilde{y})) - \ytot_{ij} \big)^2 \bigg)^{1/2}
\qquad (j\in\{1,\dots,\dimy\})
\end{align}
of the trained model's components on the totality of available data.
As our second and third quality measures for a given experimental design $\tilde{x}$, we take two different approximations of the prediction variance 
\begin{align}
V_j(x,\tilde{x}) := \Var f_j(x,\lse_f(\tilde{x},\tilde{y}(\tilde{x})))
\end{align}
of the estimated model's components based on $\tilde{x}$, namely a linearization-based approximation $V_j^{\lin}$ on the one hand and a sampling-based approximation $V_j^{\sam}$ on the other hand. 
Specifically, the linearization-based prediction variance $V_j^{\lin}(x,\tilde{x})$ 
is defined as the variance of the estimated linearized model at $x$, that is, 
\begin{align} \label{eq:linearization-based-prediction-variance}
V_j^{\lin}(x,\tilde{x}) := \Var f^{\lin}_{\ol{\theta},j}(x,\lse_{f^{\lin}_{\ol{\theta}}}(\tilde{x},\tilde{y}(\tilde{x}))) 
\qquad (j\in\{1,\dots,\dimy\}),
\end{align}
where $x \in \X$ is an arbitrary input point, $\tilde{y}(\tilde{x}):= \tilde{f}(\tilde{x},\theta^*) + \tilde{\eps}(\tilde{x})$ is the random observation model~\eqref{eq:random-observation-model}, 
and the model is linearized around $\ol{\theta} := \lse_f(\txtot, \tytot)$, which is the best parameter estimate based on the totality of available data. In view of~\eqref{eq:linearized-model} and~\eqref{eq:exact-covariance}, 
\begin{align}
V_j^{\lin}(x,\tilde{x}) = \nabla_{\theta} f_j(x,\ol{\theta})^\top M_f(\tilde{x},\ol{\theta})^{-1} \nabla_{\theta} f_j(x,\ol{\theta}).
\end{align}
Additionally, the sampling-based prediction variance $V_j^{\sam}(x,\tilde{x})$ 
is defined as the sample variance of the estimated nonlinear model at $x$, that is, 
\begin{gather} \label{eq:sampling-based-prediction-variance}
V_j^{\sam}(x,\tilde{x}) := \frac{1}{\nsam} \sum_{s=1}^{\nsam} \big( f_j(x,\lse_f(\tilde{x},\tilde{y}^{(s)})) - \mu_j^{\sam}(x,\tilde{x}) \big)^2 
\qquad (j\in\{1,\dots,\dimy\}),\\
\mu_j^{\sam}(x,\tilde{x}) := \frac{1}{\nsam} \sum_{s=1}^{\nsam} f_j(x,\lse_f(\tilde{x},\tilde{y}^{(s)})),
\end{gather}
where $x \in \X$ is an arbitrary input point, $\tilde{y}^{(s)}$ are $\nsam$ samples independently drawn from the random observation model $\tilde{y}(\tilde{x},\ol{\theta}) := \tilde{f}(\tilde{x},\ol{\theta}) + \tilde{\eps}(\tilde{x})$ with $\ol{\theta} := \lse_f(\txtot, \tytot)$. Since this $\ol{\theta}$ is the best parameter estimate based on the totality of available data, the employed observation model $\tilde{y}(\tilde{x},\ol{\theta})$ is the best available substitute for the true observation model~\eqref{eq:random-observation-model}. 
Similarly to~\eqref{eq:rmse-definition}, we will use the square roots 
\begin{align} \label{eq:linearization-and-sampling-based-prediction-uncertainties}
\sigma_j^{\lin}(x,\tilde{x}) := \big( V_j^{\lin}(x,\tilde{x}) \big)^{1/2}
\qquad \text{and} \qquad
\sigma_j^{\sam}(x,\tilde{x}) := \big( V_j^{\sam}(x,\tilde{x}) \big)^{1/2}
\end{align}
of the prediction variances~\eqref{eq:linearization-based-prediction-variance} and~\eqref{eq:sampling-based-prediction-variance} at given input points $x \in \X$. In fact, in  in our comparisons of different designs, we will actually use the $x$-independent worst-case uncertainties
\begin{align} \label{eq:linearization-and-sampling-based-worst-case-prediction-uncertainties}
\sigma_j^{\lin}(\tilde{x}) := \max_{x\in\X} \sigma_j^{\lin}(x,\tilde{x})
\qquad \text{and} \qquad
\sigma_j^{\sam}(\tilde{x}) :=  \max_{x\in\X} \sigma_j^{\sam}(x,\tilde{x}).
\end{align}
We will refer to~\eqref{eq:linearization-and-sampling-based-worst-case-prediction-uncertainties} as the worst-case linearization- and the worst-case sampling-based prediction uncertainties, respectively.
\smallskip

It is clear that for substantially nonlinear models $f$, the sampling-based prediction uncertainty $\sigma_j^{\sam}(x,\tilde{x})$ 
is a better approximation to the exact prediction uncertainty than the linearization-based prediction uncertainty $\sigma_j^{\lin}(x,\tilde{x})$, just because the former does not rely on any linearization. 
It is also clear by the very definitions~\eqref{eq:rmse-definition} and~\eqref{eq:linearization-and-sampling-based-prediction-uncertainties} of our three quality measures for experimental designs that they all depend not only on the experimental design $\tilde{x}$ alone but also on the total set of gathered experimental data $(\txtot,\tytot)$. 
It should be noticed that this dependence on $(\txtot,\tytot)$ is strongest for the root mean squared errors~\eqref{eq:rmse-definition}, while it is only indirect 
for the linearization-based prediction uncertainty~(\ref{eq:linearization-and-sampling-based-prediction-uncertainties}.a) and 
even more so for the sampling-based prediction unctertainty~(\ref{eq:linearization-and-sampling-based-prediction-uncertainties}.b). 
In order to keep the notation simple, we suppress this additional dependence on $(\txtot,\tytot)$ in the following. 

\section{Application to vapor-liquid equilibrium modeling} \label{sec:application}

In this section, we apply the proposed general sequential design method to the estimation of non-random two-liquid parameters. Specifically, we consider the binary system consisting of propanol and propyl acetate. It is known that this system exhibits narrow azeotropic behavior and therefore the estimation of its non-random two-liquid parameters is already a challenging task. With our sequential locally optimal design method, we end up with a locally optimal design consisting of $15$ experiments, four of which are repeated experiments. 
As we will see, traditional factorial design, by contrast, already requires $27$ experiments to achieve roughly the same prediction quality. 

\subsection{Considered system}

As explained above, the system considered here is the narrow azeotropic binary mixture consisting of propanol and propyl acetate. In the following, propanol will always be considered as component $1$ and propyl acetate as component $2$.
In order to estimate the binary interaction parameters of this two-component mixture, we observe it at different liquid compositions $(l,1-l)$ and at different pressures $P$ and, for all these conditions, we measure the corresponding vapor composition $(v,1-v)$ and the temperature $T$ at vapor-liquid equilibrium. In our general terminology, this means that the input quantities $x$ are the liquid mole fraction $l$ of propanol along with the system pressure $P$, while the output quantities $y$ are the vapor mole fraction $v$ of propanol along with the system temperature $T$ at vapor-liquid equilibrium. In short, 
\begin{align} \label{eq:input-and-output-quantities-case-study}
x := (l,P) \in \X := [0,1] \times [\ul{P},\ol{P}]
\qquad \text{and} \qquad 
y := (v,T) \in \Y := [0,1] \times [\ul{T},\ol{T}],
\end{align} 
where $\ul{P},\ol{P}$ and $\ul{T},\ol{T}$ determine the pressure and temperature ranges of interest. In our case, pressure and temperature are consistently measured, respectively, in units of Pascal (Pa) and Kelvin (K), and mole fractions are measured in the unit mol/mol, of course.  Also, $\ul{P} := 1 \cdot 10^5$ and $\ol{P} := 3 \cdot 10^5$ in our case. 
\smallskip

In order to model the functional relationship between the input and output quantities~\eqref{eq:input-and-output-quantities-case-study}, we assume Raoult's law in its extended version without Poynting correction. In particular, we assume non-ideal mixture behavior in the liquid phase but an ideal vapor phase (Section 13.3 of~\cite{SmVaAbSw}). 
Spelled out, this means that for every given input condition $x = (l,P) \in \X$, the equilibrium temperature $T$ is implicitly defined by the equation 
\begin{align} \label{eq:implicit-model-equation-for-T}
\frac{P_1(T)}{P} \cdot \gamma_1(l,T) l + \frac{P_2(T)}{P} \cdot  \gamma_2(l,T) (1-l) = 1,
\end{align}
and the equilibrium vapor mole fraction $v$ of propanol, in turn, is defined by the equation
\begin{align} \label{eq:implicit-model-equation-for-v}
v = \frac{P_1(T)}{P} \cdot \gamma_1(l,T) l.
\end{align}
In the above equations, $P_i(T)$ denotes the saturation vapor pressure of pure component $i$ at temperature $T$ and $\gamma_i(l,T)$ denotes the activity coefficient of component $i$ in the mixture with liquid composition $(l,1-l)$  at temperature $T$. 
Concerning the saturation vapor pressure, we assume it to be given by the Antoine model (Section~6.5 of~\cite{SmVaAbSw}), that is,
\begin{align} \label{eq:Antoine-model}
P_i(T) := 10^5 \cdot 10^{A_i - B_i/(T+C_i)}
\qquad (i \in \{1,2\}),
\end{align}
and for the Antoine parameters $A_i, B_i, C_i$ we take the estimated values from Table~\ref{tab:Antoine-parameters}.
Concerning the activity coefficients, we assume them to be given by the non-random two-liquid model of Renon and Prausnitz (Section~13.4 of~\cite{SmVaAbSw} or Section~6.15 of~\cite{PrLiAz}), 
that is, 
\begin{align} \label{eq:NRTL-model}
&\gamma_i(l,T) 
:= \exp\bigg( \frac{\sum_{k=1}^2 l_k \tau_{ki}(T) G_{ki}(T)}{\sum_{k=1}^2 l_k G_{ki}(T)} \notag\\
&\qquad \qquad \qquad \qquad \qquad \qquad + \sum_{j=1}^2 \frac{l_j G_{ij}(T)}{\sum_{k=1}^2 l_k G_{kj}(T)} \bigg( \tau_{ij}(T) - \frac{\sum_{k=1}^2 l_k \tau_{kj}(T) G_{kj}(T)}{\sum_{k=1}^2 l_k G_{kj}(T)}  \bigg) \bigg),
\end{align}
where $(l_1,l_2) := (l,1-l)$ and where $G_{ij}$ and $\tau_{ij}$ are defined by
\begin{gather}
G_{ij}(T) := \exp(-\alpha_{ij}(T)\tau_{ij}(T)) \qquad (i,j \in \{1,2\}), \label{eq:NRTL-G}\\
\tau_{ij}(T) := a_{ij} + \frac{b_{ij}}{T} \qquad (i\ne j) \qquad \text{and} \qquad \tau_{ii}(T) := 0 \qquad (i\in\{1,2\}) \label{eq:NRTL-tau}\\
\alpha_{ij}(T) := c_{ij} + d_{ij} T \qquad (i\ne j) \qquad \text{and} \qquad \alpha_{ii}(T) := 0 \qquad (i\in\{1,2\}). \label{eq:NRTL-alpha}
\end{gather}
As in the original model proposed in~\cite{RePr68}, we further assume the non-randomness terms $\alpha_{ij}$ to be symmetric and temperature-independent, that is,
\begin{align} \label{eq:NRTL-alpha-symmetric-and-T-independent}
c_{ij} = c_{ji} \qquad \text{and} \qquad d_{ij} = 0 \qquad (i,j\in\{1,2\}).
\end{align}

Since the Antoine parameters have already been estimated by a prior parameter estimation based on pure-component measurements, the parameters that remain to be estimated for our two-component mixture are the five non-random two-liquid parameters
\begin{align} \label{eq:NRTL-parameters-to-be-estimated}
\theta := (a_{12}, a_{21}, b_{12}, b_{21}, c_{12}) \in \Theta \subset \R^5.
\end{align} 
In contrast to~\cite{RePr68} and~\cite{DuAt21}, we explicitly incorporate the non-randomness parameter $c_{12} = c_{21}$ into the least-squares parameter estimation. In~\cite{RePr68}, instead, this parameter is recommended to be manually selected as one of the four values from the set $\{0.2, 0.3, 0.4, 0.47\}$, based on the characteristics of the mixture's vapor-liquid equilibrium diagram. 
\smallskip

\begin{table}[!ht]
\centering
\caption{Antoine parameters of the two components, estimated from pure-component measurements in a prior parameter estimation step.}
\label{tab:Antoine-parameters}
\begin{tabular}{lccc}
\toprule
Component &  $A_i$ &  $B_i$ & $C_i$ \\
\midrule
propanol (1)			& 4.65413     & 1292.869		& -91.992 \\
propyl acetate (2)	& 3.84871     & 1088.392 	& -90.571 \\
\bottomrule	
\end{tabular} 
\end{table}

Summarizing, the parametric model $f: \X \times \Theta \to \Y$ 
for the functional relationship between our input and output quantities~\eqref{eq:input-and-output-quantities-case-study} is defined by the equations~\eqref{eq:implicit-model-equation-for-T} and~\eqref{eq:implicit-model-equation-for-v} in conjunction with~\eqref{eq:Antoine-model}-\eqref{eq:NRTL-alpha-symmetric-and-T-independent} and, moreover, the unknown parameters $\theta$ of $f$ that are to be estimated are given by~\eqref{eq:NRTL-parameters-to-be-estimated}. Since~\eqref{eq:implicit-model-equation-for-T} and~\eqref{eq:implicit-model-equation-for-v} are algebraic equations of the form~\eqref{eq:implicit-model-definition}, the model $f$ is an implicit parametric model in the sense of our general framework (Section~\ref{sec:setting-and-terminology}). Also, in view of the strong nonlinearity of the non-random two-liquid model~\eqref{eq:NRTL-model} w.r.t.~$\theta$, it is clear that the resulting implicit model $f$ is strongly nonlinear as well. In fact, $f$ is intrinsically nonlinear to a substantial degree in the sense of~\cite{BaWa80, SeWi}.

\subsection{Apparatus and experimental setup}

In this section, we describe the measurement apparatus and process along with the initial experimental design.

\subsubsection{Apparatus and procedure}

All vapor-liquid equilibrium experiments for the considered two-component system were carried out in an ILUDEST\textsuperscript{\textregistered} measurement device, namely the FISCHER\textsuperscript{\textregistered} LABODEST\textsuperscript{\textregistered} VLE 602 apparatus. Also, the two components propanol (CAS 71-23-8) and propyl acetate (CAS 109-60-4) were taken from Geyer Chemsolute and from Sigma-Aldrich, respectively, and both had a purity of at least $0.995 \, \mathrm{mol}/\mathrm{mol}$. 
At the beginning of each experiment, the apparatus is cleaned, checked for leaks, and inerted. A mixture of the two  components is then poured into a heating zone and boiled at a constant heat rate and pressure. A Cottrell pump maintains direct contact between liquid and vapor to ensure equilibrium in a phase separation chamber. Vapor exiting the phase separation chamber is condensed, mixed with liquid, and returned to the heating zone. 
We assume vapor-liquid equilibrium to be reached when the pressure and the temperature of the system have not changed for at least $30$ minutes. The system pressure is measured by a Wika P30 sensor with 
an average absolute error of
\begin{align} \label{eq:P-measurement-error}
\sigma_P 
:= 50 \, \mathrm{Pa} 
\end{align}
and is controlled by valves connected to a nitrogen bottle and a vacuum pump. And the system temperature, in turn, is measured by a Wika Pt100 temperature sensor, located in the phase separation chamber, with an accuracy of
\begin{align} \label{eq:T-measurement-error}
\sigma_T := 0.03 \, \mathrm{K}.
\end{align} 
At equilibrium, the liquid and vapor phases are sampled after the phase separation chamber, and then their compositions are determined with an average absolute error of
\begin{align} \label{eq:mole-fraction-measurement-error}
\sigma_l, \sigma_v := 0.0015 \, \mathrm{mol/mol},
\end{align}
according to the procedure described next. So, for each experiment, there are measurement data for liquid composition $(l,1-l)$ and pressure $P$ on the one hand and for the vapor composition $(v,1-v)$ and the equilibrium temperature $T$ on the other hand.   

\subsubsection{Sample analytics}

In order to determine the compositions of the liquid and vapor phases at equilibrium, the samples taken from the phase separation chamber are first diluted with ethanol (CAS 64-17-5) at a $1:4$ sample-to-ethanol mass ratio and then measured on an Agilent gas chromatograph 6850 equipped with a thermal conductivity detector and a gas chromatography column, namely an Agilent HP-PLOT Q column of length $30 \, \mathrm{m}$,  inner diameter $0.32 \, \mathrm{mm}$, and film thickness $0.02 \, \mathrm{mm}$. The diluting ethanol was obtained from VWR Chemicals and its purity was at least $0.999 \, \mathrm{mol}/\mathrm{mol}$. The column inlet was heated to $230 \,^\circ\mathrm{C}$, helium was used as a carrier gas, and the split ratio was set to $40:1$. The oven stayed heated at $190 \,^\circ\mathrm{C}$ for $10$ minutes. And finally, two calibration curves for lower and higher concentrations of both components were used to provide concentration measurements with an average absolute error as specified in~\eqref{eq:mole-fraction-measurement-error} above.

\subsubsection{Initial design and measurement results}

\begin{figure}[!ht]
\centering
\includegraphics[width=0.6\columnwidth]{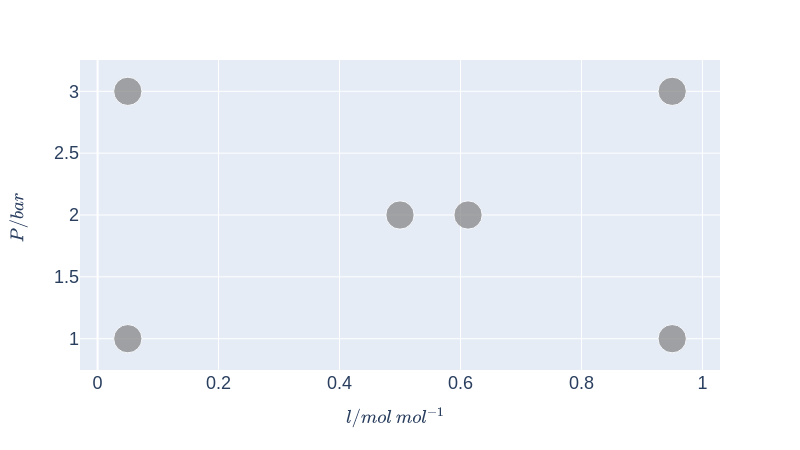}
\caption{Illustration of the experimental points from the initial design $\init$. 
}
\label{fig:initial-design}
\end{figure}

We performed the first batch of experiments according to the fractional factorial design $\init$ which consists of the four corner points and the center point of the box 
\begin{align} \label{eq:bounds-for-factorial-design}
[\ul{l}, \ol{l}] \times [\ul{P},\ol{P}] := [0.05, 0.95] \times [1\cdot 10^5, 3\cdot 10^5]
\subset \X
\end{align}
and of one additional point next to the center point (Figure~\ref{fig:initial-design}). Specifically, the initial design $\init$ consists of the six $(l,P)$ points listed in Table~\ref{tab:initial-design-and-results}. With its $6 > 5 = \dimp$ 
experimental points, the initial design of course satisfies the necessary condition~\eqref{eq:necessary-condition-for-invertibility} for the invertibility of the information matrices of $f$ at $\init$. 
It is important to notice that the experimental points of $\init$ -- and, for that matter, of any other experimental design -- cannot be realized exactly in the lab. And so there is always a mismatch between planned and actually performed experiments (which is often referred to as errors in variables in the context of parameter estimation~\cite{SeWi}). 
In particular, this is true for the planned and the actually realized liquid compositions. Apart from the planned experimental points $(l,P)$, we therefore also list the actual measurement locations $(l',P')$ in Table~\ref{tab:initial-design-and-results}. And, along with them, we list the measured values $(v,T)$ of our output quantities.   

\begin{table}[!ht]
\centering
\caption{Initial experiments conducted in the study.}
\label{tab:initial-design-and-results}
\begin{tabular}{lrrrrrrrr}
\toprule
{} &  $l $ &  $l'$ & $P$ &  $P'$ &  $v$ &     $T$ &  $\sigma_v$ &  $\sigma_T$ \\
\midrule
\init &           0.050000 &                     0.0456 &  100000.0 &      99990.0 &            0.0813 &  372.21 &                               0.0015 &                    0.03 \\
\init &           0.050000 &                     0.6961 &  300000.0 &     299970.0 &            0.7243 &  401.50 &                               0.0015 &                    0.03 \\
\init &           0.500000 &                     0.4466 &  200000.0 &     199950.0 &            0.5257 &  389.12 &                               0.0015 &                    0.03 \\
\init &           0.950000 &                     0.9728 &  100000.0 &      99990.0 &            0.9611 &  369.16 &                               0.0015 &                    0.03 \\
\init &           0.950000 &                     0.9642 &  300000.0 &     299981.6 &            0.9546 &  401.46 &                               0.0015 &                    0.03 \\
\init &           0.612500 &                     0.6426 &  200000.0 &     199950.0 &            0.6730 &  388.14 &                               0.0015 &                    0.03 \\                             
\bottomrule
\end{tabular}
\end{table}

\subsection{Application-specific settings for parameter estimation and optimal design} 

As a general modeling and optimization framework supporting automatic differentiation, we choose the software package Casadi~\cite{AnDi19}. 
Additionally, in the present application example, we use the following specific settings for the occurring parameter estimation and optimal design problems. 
In keeping with the accuracies~\eqref{eq:T-measurement-error} and~\eqref{eq:mole-fraction-measurement-error} of the employed measurement devices, we take the covariance matrix of the measurement errors to be
\begin{align} \label{eq:covariance-of-measurement-errors-case-study}
\varsigma 
:=
\begin{pmatrix}
\sigma_1^2 & 0\\
0 & \sigma_2^2
\end{pmatrix} 
=
\begin{pmatrix}
\sigma_v^2 & 0\\
0 & \sigma_T^2
\end{pmatrix} 
=
\begin{pmatrix}
2.25 \cdot 10^{-6} & 0\\
0 & 9 \cdot 10^{-4}
\end{pmatrix} 
\in \Rpd^{\dimy\times\dimy}.
\end{align} 
In particular, by the diagonality of this matrix, we implicitly assume the measurement errors for vapor composition and the measurement errors for temperature to be independent of each other. It should be recalled from~\eqref{eq:information-matrix-for-designs} and~\eqref{eq:lse-problem} that the covariance of the measurement errors enters all least-squares estimation problems and -- via the information matrices -- also all locally optimal experimental design problems. And so the above covariance matrix~\eqref{eq:covariance-of-measurement-errors-case-study} features in all optimization problems of our application example. 

\subsubsection{Settings for parameter estimation}

Since the covariance matrix~\eqref{eq:covariance-of-measurement-errors-case-study} of the measurement errors is diagonal, all least-squares estimation problems in our case study take the special form
\begin{align} \label{eq:wlse-problem-case-study}
\min_{\theta \in \Theta} \sum_{i=1}^n \sum_{j=1}^{2} (f_j(x_i,\theta)-y_{ij})^2/\sigma_j^2
\end{align}
of weighted least-squares problems. In order to solve these problems, we use the interior-point optimizer  Ipopt~\cite{WaBi05} (version 3.10.4) -- installed by Casadi -- together with the linear solver MA57~\cite{Du04, HSL}. Specifically, the optimizer settings were chosen according to the default settings of Casadi (version~3.5.5). All model derivatives required for the solution of~\eqref{eq:wlse-problem-case-study} are provided by Casadi's automatic differentiation functionality. 
In order to avoid ending up with suboptimal local minima for~\eqref{eq:wlse-problem-case-study}, we use a simple multistart procedure which randomly selects candidate start values $\theta^0$ and then filters them according to their objective function value. Specifically, if the objective function value at a candidate start value $\theta^0$ is sufficiently small, it is accepted as a start value for Ipopt; otherwise it is rejected. 

\subsubsection{Settings for optimal design}

In all locally optimal experimental design problems of our case study, we use the D-criterion as the underlying design criterion $\Psi$, that is, 
\begin{align} \label{eq:design-criterion-case-study}
\Psi := \Psi_0
\end{align}
as defined in (\ref{eq:def-D-and-A-criterion}.a). 
%
As the importance factor for already performed experiments, we consistently use $\alpha := 0.5$ throughout the case study. In this manner, we ensure that previous and prospective experiments are given equal importance throughout the entire  sequential design procedure. As the maximal tolerable optimality error for the locally $\eps$-optimal weighted designs computed by Algorithm~\ref{algo:YBT}, 
we choose $\eps := 5 \cdot 10^{-5}$. In order to solve the discretized versions~\eqref{eq:discretized-design-problem} of the weight optimization problem~\eqref{eq:2-stage-local-OED}, we use the solver SLSQP from scipy.optimize~\cite{SciPy} 
with its default settings (version 12.1.0). (It should be noted that other solvers could have been used just as well. After all, as has been explained after Algorithm~\ref{algo:YBT}, the discretized problems~\eqref{eq:discretized-design-problem} are just convex weight optimization problems with only as many variables as the size of the discretization $\ul{X}^k$.) All model derivatives required for the computation of the information matrices in~\eqref{eq:discretized-design-problem} and~\eqref{eq:violator-problem} are provided by Casadi's automatic differentiation functionality. 
%
As the minimal total weight of the candidate experiments that remain after our sorting-out procedure (Algorithm~\ref{algo:2-stage-locally-optimal-design}), we choose $\ul{w}^+ := 0.95$. And as the maximal number of new experiments per iteration of our sequential design algorithm, 
we choose $\ol{n}^+ := 3$. We make this choice because, from a comparison perspective, it fits well with our reference factorial designs of the sizes $3 \cdot 3$, $5 \cdot 3$, and $9 \cdot 3$, respectively. 
%
In order to initialize our sequential design procedure (Algorithm~\ref{algo:sequential-optimal-design}), we choose the $6$-point fractional factorial design 
\begin{align} \label{eq:seq-algo-inputs-1-case-study}
\tilde{x}^0 := \init
\end{align}
introduced above (Figure~\ref{fig:initial-design} and Table~\ref{tab:initial-design-and-results}). And finally, the maximal number of iterations and the progress tolerance of our sequential design algorithm are chosen to be $\ol{n} := 27$ and $\delta := 0.1$ in our case study. 
Summarizing, we apply Algorithm~\ref{algo:sequential-optimal-design} with the following inputs:
\begin{gather} 
\alpha := 0.5, \qquad \eps := 0.5 \cdot 10^{-4}, \qquad \ul{w}^+ := 0.95, \label{eq:seq-algo-inputs-2-case-study}\\
\ol{n}^+ := 3, \qquad \ol{n} := 27, \qquad \delta := 0.1. \label{eq:seq-algo-inputs-3-case-study}
\end{gather}

\subsection{Considered designs, measurement results, and paremeter estimates}

After performing the six experiments from the initial fractional factorial design $\init$ (Table~\ref{tab:initial-design-and-results}), we continue to measure the output quantities~\eqref{eq:input-and-output-quantities-case-study}  
at a number of additional input points, namely first at a sequence
\begin{align} \label{eq:factorial-designs}
\fed^1 = \fed^0 \,\&\, \fed^{0,+}, \qquad \fed^2 = \fed^1 \,\&\, \fed^{1,+}, \qquad \fed^3 = \fed^2 \,\&\, \fed^{2,+}
\end{align} 
of factorial designs $\fed^1, \fed^2, \fed^3$, and then at a sequence
\begin{align} \label{eq:optimal-designs}
\oed^1 = \oed^0 \,\&\, \oed^{0,+}, \qquad \oed^2 = \oed^1 \,\&\, \oed^{1,+}, \qquad \oed^3 = \oed^2 \,\&\, \oed^{2,+}
\end{align}
of locally optimal designs $\oed^1, \oed^2, \oed^3$. We used the notation~\eqref{eq:combined-design} for refined designs here. In particular, all designs from~\eqref{eq:factorial-designs} arise by a  successive refinement of the $5$-point fractional factorial design $\fed^0$ consisting of the first five experiments from Table~\ref{tab:initial-design-and-results} (all points in Figure~\ref{fig:initial-design} except for the off-center point). Similarly, all designs from~\eqref{eq:optimal-designs} arise by a  successive refinement of the $6$-point fractional factorial design $\oed^0 = \init$ consisting of the six experiments from Table~\ref{tab:initial-design-and-results}.
\smallskip

\begin{figure}
\centering
\includegraphics[width=0.6\columnwidth]{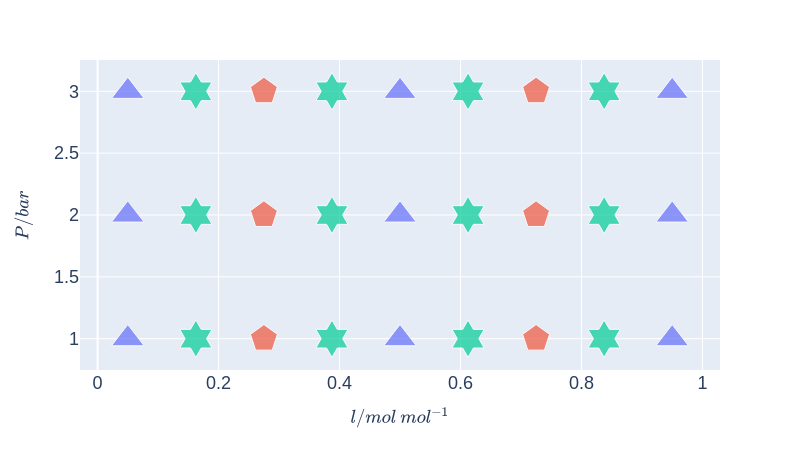}
\caption{Illustration of the experimental points from the factorial designs $\fed^1$ (blue triangles), $\fed^2 = \fed^1 \,\&\, \fed^{1,+}$, and $\fed^3 = \fed^2 \,\&\, \fed^{2,+}$. The experimental points of the designs $\fed^{1,+}$ and $\fed^{2,+}$ are the red pentagons and the green stars, respectively.
}
\label{fig:factorial-designs}
\end{figure}

As the designs $\fed^1, \fed^2, \fed^3$ in~\eqref{eq:factorial-designs}, we choose full factorial designs on the $9 \cdot 3$ equidistant grid
\begin{align} \label{eq:design-space-FED}
\ul{\X}^{\mathrm{fed}} := \big\{ (\ul{l} + i (\ol{l}-\ul{l})/8, \ul{P} + j(\ol{P}-\ul{P})/2): i \in \{0,1,\dots,8\} \text{ and } j \in \{0,1,2\} \big\} \subset \X,
\end{align}
where $\ul{l} := 0.05$ and $\ol{l} := 0.95$ are the mole fraction bounds from~\eqref{eq:bounds-for-factorial-design} 
and $\ul{P}, \ol{P}$ are the pressure bounds from~\eqref{eq:input-and-output-quantities-case-study}.  
Specifically, $\fed^1, \fed^2, \fed^3$ are the $3 \cdot 3$, $5 \cdot 3$, and $9 \cdot 3$ full factorial designs depicted in Figure~\ref{fig:factorial-designs} and defined by the $l$- and $P$-columns in Table~\ref{tab:factorial-designs-and-results}.
\smallskip

\begin{table}
\centering
\caption{Sequence of factorial experiments~\eqref{eq:factorial-designs} conducted in the study (following the five  initial experiments of $\fed^0$, that is, the first five experiments from $\init$).}
\label{tab:factorial-designs-and-results}
\begin{tabular}{lrrrrrrrr}
\toprule
{} &  $l $ &  $l'$ & $P$ &  $P'$ &  $v$ &     $T$ &  $\sigma_v$ &  $\sigma_T$ \\
\midrule
$\fed^{0,+}$            &           0.050000 &                     0.0444 &  200000.0 &     199950.0 &            0.0815 &  396.06 &                               0.0015 &                    0.03 \\
$\fed^{0,+}$           &           0.500000 &                     0.4469 &  100000.0 &      99990.0 &            0.4979 &  367.32 &                               0.0015 &                    0.03 \\
$\fed^{0,+}$             &           0.500000 &                     0.3538 &  300000.0 &     299981.6 &            0.4560 &  404.32 &                               0.0015 &                    0.03 \\
$\fed^{0,+}$            &           0.950000 &                     0.9685 &  200000.0 &     199950.0 &            0.9525 &  388.73 &                               0.0015 &                    0.03 \\
$\fed^{1,+}$            &           0.275000 &                     0.2440 &  100000.0 &      99990.0 &            0.3341 &  368.70 &                               0.0015 &                    0.03 \\
$\fed^{1,+}$            &           0.275000 &                     0.2429 &  200000.0 &     199950.0 &            0.3431 &  391.58 &                               0.0015 &                    0.03 \\
$\fed^{1,+}$            &           0.275000 &                     0.1168 &  300000.0 &     299981.6 &            0.1933 &  409.82 &                               0.0015 &                    0.03 \\
$\fed^{1,+}$           &           0.725000 &                     0.6974 &  100000.0 &      99990.0 &            0.6750 &  367.12 &                               0.0015 &                    0.03 \\
$\fed^{1,+}$           &           0.725000 &                     0.6958 &  200000.0 &     199950.0 &            0.7099 &  387.82 &                               0.0015 &                    0.03 \\
$\fed^{1,+}$           &           0.725000 &                     0.5809 &  300000.0 &     299981.6 &            0.6499 &  402.05 &                               0.0015 &                    0.03 \\
$\fed^{2,+}$            &           0.162500 &                     0.1503 &  100000.0 &      99990.0 &            0.2385 &  369.95 &                               0.0015 &                    0.03 \\
$\fed^{2,+}$            &           0.162500 &                     0.1486 &  200000.0 &     199950.0 &            0.2432 &  393.22 &                               0.0015 &                    0.03 \\
$\fed^{2,+}$           &           0.162500 &                     0.0454 &  300000.0 &     299981.6 &            0.0815 &  411.82 &                               0.0015 &                    0.03 \\
$\fed^{2,+}$           &           0.387500 &                     0.3707 &  100000.0 &      99990.0 &            0.4378 &  367.46 &                               0.0015 &                    0.03 \\
$\fed^{2,+}$           &           0.387500 &                     0.3697 &  200000.0 &     199950.0 &            0.4430 &  389.46 &                               0.0015 &                    0.03 \\
$\fed^{2,+}$           &           0.387500 &                     0.2346 &  300000.0 &     299981.6 &            0.3424 &  406.86 &                               0.0015 &                    0.03 \\
$\fed^{2,+}$           &           0.612500 &                     0.6037 &  100000.0 &      99990.0 &            0.6026 &  367.05 &                               0.0015 &                    0.03 \\
$\fed^{2,+}$		 &           0.612500 &                     0.6426 &  200000.0 &     199950.0 &            0.6730 &  388.14 &                               0.0015 &                    0.03 \\ 
$\fed^{2,+}$          &           0.612500 &                     0.4356 &  300000.0 &     299981.6 &            0.5353 &  403.59 &                               0.0015 &                    0.03 \\
$\fed^{2,+}$            &           0.837500 &                     0.9168 &  100000.0 &      99990.0 &            0.8869 &  368.41 &                               0.0015 &                    0.03 \\
$\fed^{2,+}$            &           0.837500 &                     0.9037 &  200000.0 &     199950.0 &            0.8887 &  388.39 &                               0.0015 &                    0.03 \\
$\fed^{2,+}$         &           0.837500 &                     0.9067 &  300000.0 &     299981.6 &            0.8903 &  401.36 &                               0.0015 &                    0.03 \\
\bottomrule
\end{tabular}
\end{table}

\begin{figure}
\centering
\includegraphics[width=0.6\columnwidth]{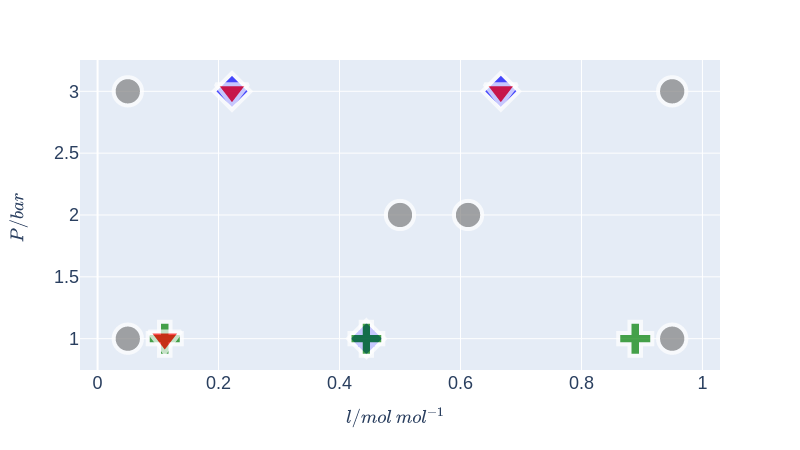}
\caption{Illustration of the experimental points from the optimal designs $\oed^0 = \init$ (gray circles), and $\oed^1 = \oed^0 \,\&\, \oed^{0,+}$, $\oed^2 = \oed^1 \,\&\, \oed^{1,+}$, and $\oed^3 = \oed^2 \,\&\, \oed^{2,+}$. The experimental points of the designs $\oed^{0,+}$, $\oed^{1,+}$ and $\oed^{2,+}$ are the blue diamonds, the green pluses and the red triangles, respectively.
}
\label{fig:optimal-designs}
\end{figure}

As the designs $\oed^1, \oed^2, \oed^3$ in~\eqref{eq:optimal-designs}, we choose the locally optimal designs computed in the course of our sequential design algorithm (Algorithm~\ref{algo:sequential-optimal-design}) with design criterion~\eqref{eq:design-criterion-case-study}, inputs~\eqref{eq:seq-algo-inputs-1-case-study}-\eqref{eq:seq-algo-inputs-3-case-study}, and with the $10 \cdot 10$ equidistant grid
\begin{align} \label{eq:design-space-OED}
\ul{\X}^{\mathrm{oed}} := \big\{ (i/9, \ul{P} + j(\ol{P}-\ul{P})/9): i, j \in \{0,1,\dots,9\} \big\} \subset \X,
\end{align}
as the underlying design space, where again $\ul{P}, \ol{P}$ are the pressure bounds from~\eqref{eq:input-and-output-quantities-case-study}. With these settings, the sequential design algorithm terminates after just three iterations because all experimental points of the design $\oed^{2,+}$ proposed in the third iteration already belong to the previous design $\oed^2$. And therefore we end up with a sequence of just three locally optimal designs $\oed^1, \oed^2, \oed^3$. 
See Figure~\ref{fig:optimal-designs} and the $l$- and $P$-columns in Table~\ref{tab:optimal-designs-and-results} for the experimental points making up these designs. In particular, the designs $\oed^1, \oed^2, \oed^3$ have  the sizes $9$, $12$ and $15$. Also, while $\oed^1$ is a design without replications, $\oed^2$ and $\oed^3$ contain replicated -- or more precisely, duplicated -- experiments. Specifically, $\oed^1$ contains one duplicated experiment and $\oed^2$ contains four duplicated experiments. See the overlaps in the figure. 
\smallskip

\begin{table}
\centering
\caption{Sequence of optimal experiments~\eqref{eq:optimal-designs} conducted in the study (following the six  initial experiments from $\oed^0 = \init$).}
\label{tab:optimal-designs-and-results}
\begin{tabular}{lrrrrrrrr}
\toprule
{} &  $l $ &  $l'$ & $P$ &  $P'$ &  $v$ &     $T$ &  $\sigma_v$ &  $\sigma_T$ \\
\midrule
$\oed^{0,+}$             &           0.222222 &                     0.1731 & 300000.0   &  299981.6    &            0.2748 &  407.73 &                               0.0015 &                    0.03 \\
$\oed^{0,+}$              &           0.444444 &                     0.3878 & 100000.0   & 99990.0     &            0.4561 &  367.49 &                               0.0015 &                    0.03 \\
$\oed^{0,+}$             &           0.666667 &                     0.6520 & 300000.0  & 299981.6     &            0.6952 &  401.62 &                               0.0015 &                    0.03 \\
$\oed^{1,+}$             &           0.111111 &                     0.1018 &   100000.0 &      99990.0 &            0.1721 &  370.81 &                               0.0015 &                    0.03 \\
$\oed^{1,+}$             &           0.444444 &                     0.4518 &   100000.0 &      99990.0 &            0.4992 &  367.23 &                               0.0015 &                    0.03 \\
$\oed^{1,+}$             &           0.888889 &                     0.8845 &   100000.0 &      99990.0 &            0.8459 &  368.03 &                               0.0015 &                    0.03 \\
$\oed^{2,+}$              &           0.111111 &                     0.1094 &   100000.0 &      99990.0 &            0.1803 &  370.68 &                               0.0015 &                    0.03 \\
$\oed^{2,+}$             &           0.222222 &                     0.2385 &  300000.0 &     300000.0 &            0.3466 &  406.26 &                               0.0015 &                    0.03 \\
$\oed^{2,+}$            &           0.666667 &                     0.7372 &  300000.0 &     300000.0 &            0.7586 &  401.28 &                               0.0015 &                    0.03 \\
\bottomrule
\end{tabular}
\end{table}

Just like for the initial design, Tables~\ref{tab:factorial-designs-and-results} and~\ref{tab:optimal-designs-and-results} not only specify the planned experimental designs $\fed^k$ and $\oed^k$, but also detail the input conditions $(l',P')$ that were actually realized in the lab experiments, along with the corresponding measurement results $(v,T)$ for the outputs. In particular, we see that the planned replications in $\oed^2$ and in $\oed^3$ represented by overlaps in Figure~\ref{fig:optimal-designs} are fanned out a bit, in reality. 
\smallskip

\begin{figure}[!ht]
\centering
\includegraphics[width=0.6\columnwidth]{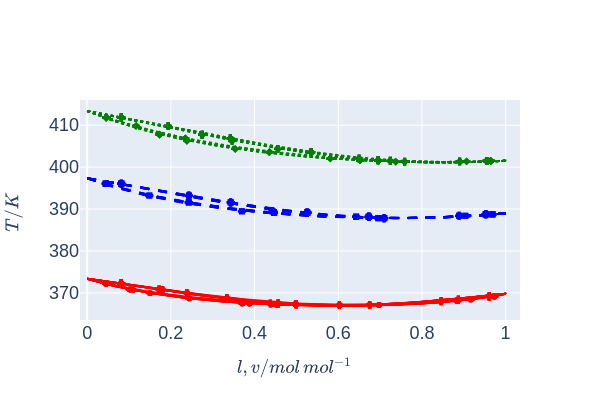}
\caption{The predicted and measured bubble- and dew-point temperatures for propanol and propyl acetate at the pressures $1\, \mathrm{bar}$ (bottom-most curves), $2\, \mathrm{bar}$ (middle curves), and $3\, \mathrm{bar}$ (top-most curves). The predictions are based on the final least-squares estimate $\lse^\tot$ from Table~\ref{tab:lse-tot}.}
\label{fig:VLE-tot}
\end{figure}

In the following, we write $(\txtot,\tytot)$ for the totality of the collected experimental data, that is, for the combination of all the $(l,P)$- and $(v,T)$-values from Tables~\ref{tab:initial-design-and-results}, 
\ref{tab:factorial-designs-and-results} and~\ref{tab:optimal-designs-and-results}.
We also write $\tot := \txtot$ for the total combined design. In formulas, 
\begin{align} \label{eq:total-design}
\txtot := \tot := \fed^0 \,\&\, \fed^{0,+} \,\&\, \fed^{1,+} \,\&\, \fed^{2,+} \,\&\, \oed^{0,+} \,\&\, \oed^{1,+} \,\&\, \oed^{2,+}, 
\end{align}
so that the total combined design $\tot$ has the size $\ntot = 36$. Table~\ref{tab:lse-tot} records our final least-squares estimate for the non-random two liquid parameters of propanol and propyl acetate, that is, the least-squares estimate $\lse^\tot := \lse_f(\txtot,\tytot)$ based on all collected data $(\txtot,\tytot)$. Figure~\ref{fig:VLE-tot} shows the corresponding predicted bubble- and dew-point temperature curves for the three pressures considered in the experiments, along with the measured temperature values. As can be seen, the predictions are in perfect agreement with the measurements. Additionally, the figure nicely illustrates the narrow azeotropic behavior of our mixture.

\begin{table}[!ht]
\centering
\caption{The final least-squares estimate $\lse^\tot := \lse_f(\txtot,\tytot)$ for the non-random two liquid parameters of propanol (1) and propyl acetate (2), based on all collected data $(\txtot,\tytot)$.}
\label{tab:lse-tot}
\begin{tabular}{ccccc}
\toprule
$a_{1,2}$ &  $a_{2,1}$  & $b_{1,2}$ &  $b_{2,1}$ &  $c_{1,2} = c_{2,1}$  \\
\midrule
9.396525 &       -10.305843 &      -786.446701 &      1510.352034 &         0.010000	 \\
\bottomrule
\end{tabular} 
\end{table}

\subsection{Comparison of the considered designs and of the resulting models}

In this section, we finally compare the prediction quality of the models trained on the factorial designs $\fed^k$ with the prediction quality of the models trained on the locally optimal designs $\oed^k$ computed by our sequential design algorithm. In order to do so, we use the three prediction quality measures~\eqref{eq:rmse-definition}, (\ref{eq:linearization-and-sampling-based-worst-case-prediction-uncertainties}.a) and~(\ref{eq:linearization-and-sampling-based-worst-case-prediction-uncertainties}.b) from Section~\ref{sec:quality-measures}. As has been pointed out there, these prediction quality measures all depend on the totality $(\txtot,\tytot)$ of the collected experimental data as reference data, see~\eqref{eq:total-design} above. 
\smallskip

\begin{figure}[!ht]
\centering
\includegraphics[width=\columnwidth]{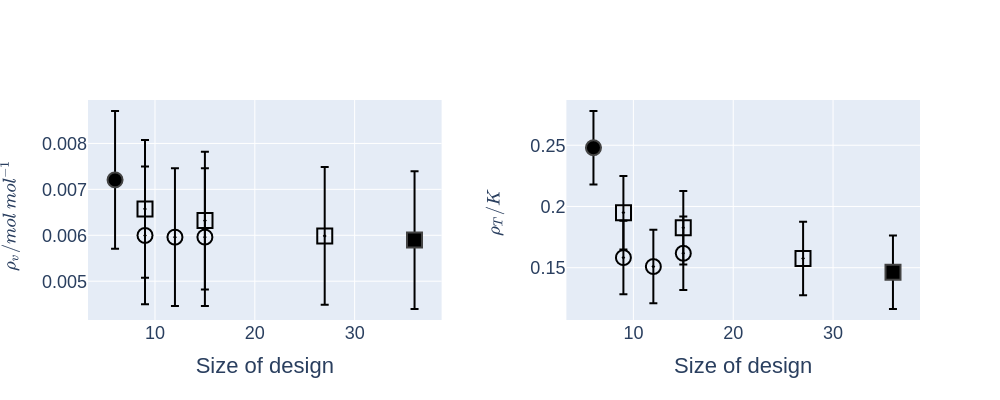}
\caption{Root mean squared prediction errors for the models trained on the different designs. The open squares and open circles represent respectively the designs $\fed^k$ and $\oed^k$, while the closed circle and the closed square represent the design $\init$ and $\tot$, respectively. 
As can be seen, $\oed^1$ with its $9$ experiments already achieves the same prediction error as $\fed^3$ with its $27$ experiments.}
\label{fig:rmse}
\end{figure}

We begin by considering the root mean squared prediction errors $\rho_j(\tilde{x})$ for the models trained on the different designs, see~\eqref{eq:rmse-definition}. Since our model has two outputs~\eqref{eq:input-and-output-quantities-case-study}, there are two such prediction errors for each design $\tilde{x}$ and we denote them by $\rho_v(\tilde{x}) := \rho_1(\tilde{x})$ and $\rho_T(\tilde{x}) := \rho_2(\tilde{x})$, respectively. Table~\ref{tab:rmse} lists these values for all considered designs $\tilde{x}$, and Figure~\ref{fig:rmse} plots them as functions of the design size.
\smallskip

\begin{table}[!ht]
\centering
\caption{Root mean squared prediction errors for the models trained on the different designs}
\label{tab:rmse}
\begin{tabular}{lrrrrrrrr}
\toprule
{} &  $\init$ &  $\oed^1$  & $\fed^1$ &  $\oed^2$ &  $\fed^2$ &     $\oed^3$ &  $\fed^3$ &  $\tot$ \\
\midrule
$\operatorname{size}(\tilde{x})$& 6      & 9		& 9		 & 12	   & 15	    & 15		 & 27	 	& 36 \\
$\rho_v(\tilde{x}) / (10^{-4})$ & 72.07	& 59.96	& 65.74	 & 59.59   & 63.21   & 59.61  & 59.86    & 58.95 \\
$\rho_T(\tilde{x}) / (10^{-2})$ & 24.80  & 15.83 & 19.50   & 15.10   & 18.26  & 16.18  & 15.75    & 14.63 \\
\bottomrule
\end{tabular} 
\end{table}

As can be seen from this figure, the prediction error $\rho_j(\oed^3)$ based on the size-$15$ locally optimal design $\oed^3$ 
is practically the same as the prediction error $\rho_j(\fed^3)$ based on the size-$27$ full factorial design $\fed^3$. 
In short, 
\begin{align} \label{eq:equal-rmse}
\rho_j(\oed^3) \approx \rho_j(\fed^3)
\qquad (j \in \{1,2\})
\end{align}
even though $\oed^3$ is only about half the size of $\fed^3$. Spelled out, \eqref{eq:equal-rmse} means that on the totality of collected data, the model trained just on the $15$ optimally designed experiments from $\oed^3$ makes essentially the same prediction error as the model trained on the $27$  factorially designed experiments from $\fed^3$.  
\smallskip

Additionally, the figure and the table show that the prediction errors $\rho_j(\oed^3)$ and $\rho_j(\fed^3)$ from~\eqref{eq:equal-rmse} are not only essentially equal to each other, but also almost equal to $\rho_j(\tot)$ which, in turn, is 
the smallest possible prediction error value on our total data set. After all, 
\begin{align}
\sum_{j=1}^2 \rho_j(\tot)^2/\sigma_j^2 
= \min_{\theta \in \Theta} \frac{1}{\ntot} \sum_{i=1}^{\ntot} \sum_{j=1}^2 (f_j(\xtot_i,\theta)-\ytot_{ij})^2/\sigma_j^2
\end{align}
by the very definition of the least-squares estimate $\lse_f(\txtot,\tytot)$ for $(\txtot,\tytot)$. 
\smallskip

And finally, the figure and the table reveal that the prediction error is essentially constant along the designs $\oed^k$ computed by our algorithm, 
while it is strictly decreasing along our factorial designs $\fed^k$. 
In particular, the terminal prediction error $\rho_j(\oed^3) \approx \rho_j(\fed^3)$ is reached already with $\oed^1$ and $\oed^2$, while it is missed with $\fed^1$ or $\fed^2$. So, from a pure prediction-error perspective, it would have made sense to terminate our algorithm even earlier. 
Yet, from a prediction-uncertainty perspective, a premature termination turns out to be not advisable, as we will see in a moment.  
\smallskip 

We continue by considering the worst-case prediction uncertainties $\sigma_j^{\lin}(\tilde{x})$ and $\sigma_j^{\sam}(\tilde{x})$ for the models trained on the different designs, see~\eqref{eq:linearization-and-sampling-based-worst-case-prediction-uncertainties}. In analogy to the root mean squared prediction errors, we denote them by $\sigma^{\lin}_v(\tilde{x}) := \sigma_1^{\lin}(\tilde{x})$, $\sigma^{\lin}_T(\tilde{x}) := \sigma_2^{\lin}(\tilde{x})$ and by $\sigma^{\sam}_v(\tilde{x}) := \sigma_1^{\sam}(\tilde{x})$, $\sigma^{\sam}_T(\tilde{x}) := \sigma_2^{\sam}(\tilde{x})$, respectively, and for the sampling-based uncertainties we used samples of the size 
\begin{align}
\nsam := 1000.
\end{align}
Tables~\ref{tab:sigma-lin} and~\ref{tab:sigma-sam} record these values for all designs $\tilde{x}$ considered here. 
\smallskip

\begin{table}[!ht]
\centering
\caption{Worst-case prediction uncertainties (linearization-based) for the models trained on the different designs}
\label{tab:sigma-lin}
\begin{tabular}{lrrrrrrrr}
\toprule
{} &  $\init$ &  $\oed^1$  & $\fed^1$ &  $\oed^2$ &  $\fed^2$ &     $\oed^3$ &  $\fed^3$ &  $\tot$ \\
\midrule
$\operatorname{size}(\tilde{x})$ 		& 6     & 9		& 9		 & 12	 & 15	& 15		& 27	 	& 36 \\
$\sigma^{\lin}_v(\tilde{x}) / (10^{-4})$ & 67.90	& 33.08	& 32.08	 & 28.11  & 27.46  & 25.47  & 24.92  & 23.07 \\
$\sigma^{\lin}_T(\tilde{x}) / (10^{-2})$ & 34.14  & 10.79  & 10.87    & 10.40  & 8.96  & 8.26  & 8.37  & 7.85 \\
\bottomrule
\end{tabular} 
\end{table}

\begin{table}[!ht]
\centering
\caption{Worst-case prediction uncertainties (sampling-based) for the models trained on the different designs}
\label{tab:sigma-sam}
\begin{tabular}{lrrrrrrrr}
\toprule
{} &  $\init$ &  $\oed^1$ & $\fed^1$ &  $\oed^2$ &  $\fed^2$ &     $\oed^3$ &  $\fed^3$ &  $\tot$ \\
\midrule
$\operatorname{size}(\tilde{x})$ 		& 6     & 9		& 9		 & 12	 & 15	& 15		& 27	 	& 36 \\
$\sigma^{\sam}_v(\tilde{x}) / (10^{-4})$ & 17.3	& 10.78	& 9.84	 & 6.48  & 6.98  & 5.60  & 4.49  & 3.71 \\
$\sigma^{\sam}_T(\tilde{x}) / (10^{-2})$ & 7.56  & 3.08  & 3.75    & 2.53  & 2.39  & 1.94  & 1.36  & 1.13 \\
\bottomrule
\end{tabular} 
\end{table}

As in the prediction-error case considered above, it can be seen from the prediction-uncertainty tables that the worst-case prediction uncertainties $\sigma_j^{\lin}(\oed^3)$ and $\sigma_j^{\sam}(\oed^3)$ based on the size-$15$ locally optimal design $\oed^3$ are practically the same as the worst-case prediction uncertainties $\sigma_j^{\lin}(\fed^3)$ and $\sigma_j^{\sam}(\fed^3)$ based on the size-$27$ full factorial design $\fed^3$. 
In short, 
\begin{align} \label{eq:equal-prediction-uncertainties}
\sigma_j^{\lin}(\oed^3) \approx \sigma_j^{\lin}(\fed^3)
\qquad \text{and} \qquad
\sigma_j^{\sam}(\oed^3) \approx \sigma_j^{\sam}(\fed^3)
\end{align}
despite the substantially different sizes of $\oed^3$ and $\fed^3$. In words, 
\eqref{eq:equal-prediction-uncertainties} means that the predictions of the model trained just on the $15$ optimally designed experiments from $\oed^3$ can be trusted basically as much as the predictions of the model trained on the $27$ factorially designed experiments from $\fed^3$.
\smallskip

As opposed to the prediction-error case, the terminal prediction uncertainties from~\eqref{eq:equal-prediction-uncertainties} are no longer reached already for $\oed^1$ or $\oed^2$. Instead, the prediction uncertainties are strictly improving not only along our factorial designs $\fed^k$, but also along the optimal designs $\oed^k$ computed by our algorithm. And therefore, from a prediction-uncertainty perspective, a premature termination of our algorithm would be disadvantageous.
\smallskip 

\begin{figure}[!ht]
\centering
\includegraphics[width=\columnwidth]{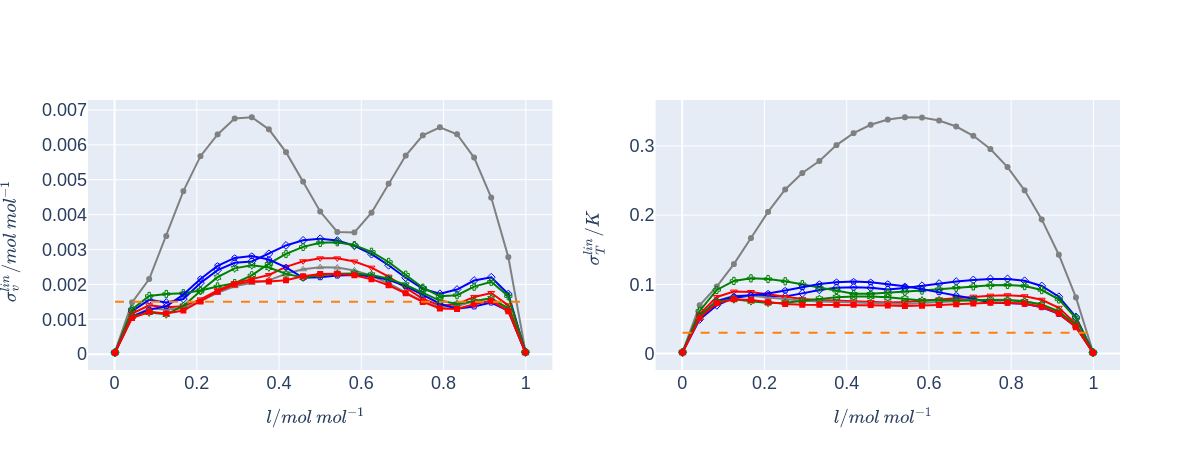}
\caption{Worst-pressure-case prediction uncertainties (linearization-based) for the models trained on the different designs. The prediction uncertaintiy curves~\eqref{eq:worst-pressure-case-prediction-uncertainties} for $\init$, $\oed^1$, $\fed^1$, $\oed^2$, $\fed^2$, $\oed^3$, $\fed^3$ and $\tot$ are represented, respectively, by the gray circles, the blue diamonds, the green pluses, the blue circles, the red triangles, the green stars, the gray triangles, and the red squares. The orange lines indicate the measurment uncertainties $\sigma_v$ and $\sigma_T$. Comparing the maxima of the green-stars and the gray-triangles curves, we see that the worst-case prediction uncertainty achieved with the size-$15$ design $\oed^3$ is practically the same as the worst-case prediction uncertainty achieved with the size-$27$ design $\fed^3$.}
\label{fig:sigma-lin}
\end{figure}

\begin{figure}[!ht]
\centering
\includegraphics[width=\columnwidth]{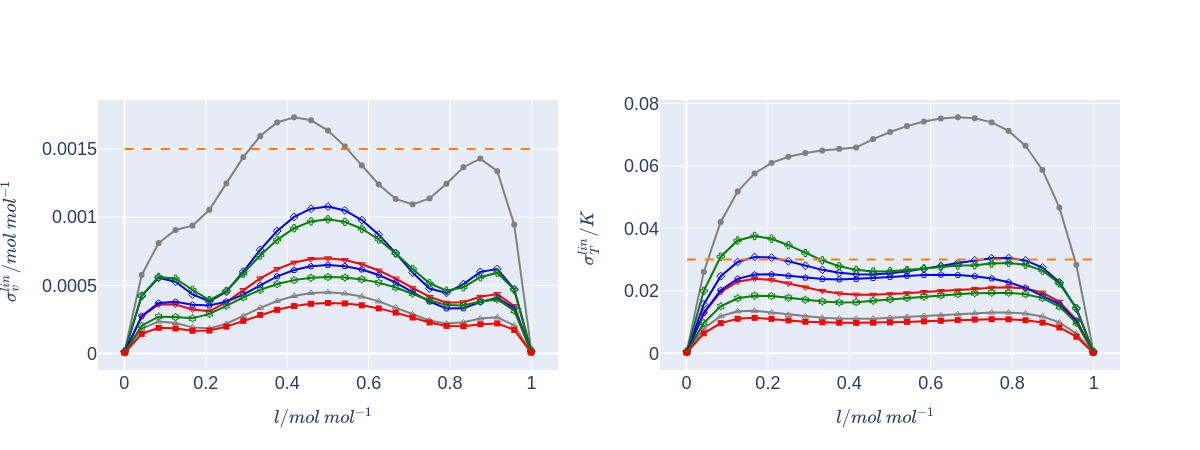}
\caption{Worst-pressure-case prediction uncertainties (sampling-based) for the models trained on the different designs. The prediction uncertainty curves~\eqref{eq:worst-pressure-case-prediction-uncertainties} for $\init$, $\oed^1$, $\fed^1$, $\oed^2$, $\fed^2$, $\oed^3$, $\fed^3$ and $\tot$ are represented, respectively, by the gray circles, the blue diamonds, the green pluses, the blue circles, the red triangles, the green stars, the gray triangles, and the red squares. The orange lines indicate the measurment uncertainties $\sigma_v$ and $\sigma_T$. Comparing the maxima of the green-stars and the gray-triangles curves, we see that the worst-case prediction uncertainty achieved with the size-$15$ design $\oed^3$ is practically the same as the worst-case prediction uncertainty achieved with the size-$27$ design $\fed^3$.}
\label{fig:sigma-sam}
\end{figure}

In order to highlight the variability of the prediction uncertainties~\eqref{eq:linearization-and-sampling-based-prediction-uncertainties}, we plot them as functions of one of our inputs $x = (l,P)$, namely the liquid mole fraction $l$ of propanol (Figures~\ref{fig:sigma-lin} and~\ref{fig:sigma-sam}). Specifically, 
the curves in Figures~\ref{fig:sigma-lin} and~\ref{fig:sigma-sam} depict the worst-pressure-case prediction uncertainties
\begin{align} \label{eq:worst-pressure-case-prediction-uncertainties}
l \mapsto \max_{P\in[\ul{P},\ol{P}]} \sigma_j^{\lin}((l,P),\tilde{x})
\qquad \text{and} \qquad
l \mapsto \max_{P\in[\ul{P},\ol{P}]} \sigma_j^{\sam}((l,P),\tilde{x})
\end{align}
as functions of $l$ for all the considered designs $\tilde{x}$. Clearly, the maxima of these worst-pressure-case curves are nothing but the worst-case prediction uncertainties $\sigma_j^{\lin}(\tilde{x})$ and $\sigma_j^{\sam}(\tilde{x})$ that we considered above (Tables~\ref{tab:sigma-lin} and~\ref{tab:sigma-sam}). 
\smallskip

We immediately observe from the figures that the worst-pressure-case prediction uncertainty curves~\eqref{eq:worst-pressure-case-prediction-uncertainties} are strongly varying with $l$, in general. And so, the prediction uncertainties~\eqref{eq:linearization-and-sampling-based-prediction-uncertainties} 
are strongly dependent on the inputs $x$ as well, in general. 
We also observe from the figures that the variability of the depicted curves~\eqref{eq:worst-pressure-case-prediction-uncertainties} becomes less and less pronounced, the larger the size of the design $\tilde{x}$. In other words, the curves~\eqref{eq:worst-pressure-case-prediction-uncertainties} get closer and closer to being constants equal to the worst-case prediction uncertainties $\sigma_j^{\lin}(\tilde{x})$ and $\sigma_j^{\sam}(\tilde{x})$. 
It should be noticed, however, that the curves~\eqref{eq:worst-pressure-case-prediction-uncertainties} will alway retain some variability -- they can never become exactly constant just because at the boundaries $l = 0$ (pure propyl acetate) and $l = 1$ (pure propanol) they take exactly the value $0$. In short,
\begin{align} \label{eq:worst-pressure-case-prediction-uncertainties-zero-at-boundaries}
\sigma_j^{\lin}((0,P),\tilde{x}) = 0 = \sigma_j^{\lin}((1,P),\tilde{x})
\qquad \text{and} \qquad 
\sigma_j^{\sam}((0,P),\tilde{x}) = 0 = \sigma_j^{\sam}((1,P),\tilde{x})
\end{align}
for every arbitrary pressure value $P$ and for every arbitrary experimental design $\tilde{x}$. Indeed, the simple reason behind~\eqref{eq:worst-pressure-case-prediction-uncertainties-zero-at-boundaries} is that for the pure-component cases $l=0$ and $l=1$, the non-random two-liquid model~\eqref{eq:NRTL-model} does not depend on the parameters $\theta$ at all. And therefore, the same is true for the solutions $T(l,P,\theta)$ and $v(l,P,\theta)$ of~\eqref{eq:implicit-model-equation-for-T} and~\eqref{eq:implicit-model-equation-for-v}. So, for the boundary values $l=0$ and $l=1$,  the model predictions
\begin{align}
f(l,P,\theta) := (v(l,P,\theta), T(l,P,\theta))^\top
\end{align}
are actually completely independent of $\theta$ and this, in conjunction with the definitions~\eqref{eq:linearization-and-sampling-based-prediction-uncertainties}, immediately implies the relations claimed in~\eqref{eq:worst-pressure-case-prediction-uncertainties-zero-at-boundaries} above and evidenced by the figures.  

\section{Conclusions and outlook}

In this paper, we have proposed a general method of sequential locally optimal experimental design for nonlinear models. It can be applied to 
estimate the parameters of nonlinear models as they abound in vapor-liquid equilibrium modeling, for instance. And these nonlinear models can be defined either by explicit or by implicit algebraic equations. 
In every iteration of our method, a whole batch of new experiments is computed in a two-stage locally optimal manner. In other words, a two-stage locally optimal experimental design problem is solved in every iteration of our method, based on the previously collected experimental data and on the corresponding least-squares estimate. As the underlying design criterion, any of the standard differentiable design criteria can be chosen, in particular, the A- and the D-criterion. Also, the number of new experiments computed in each iteration can be freely specified by the user in such a way that it fits the concrete experimental circumstances best. 
\smallskip

We have demonstrated the benefits of the general method by applying it to the estimation of non-random two-liquid parameters of a narrow azeotropic system. In particular, we have found that the proposed sequential locally optimal experimental design method requires substantially fewer experiments than traditional factorial experimental design to achieve the same model precision and prediction quality. Specifically, compared to factorial design, our sequential design method required only about half the experiments and thus helps to considerably reduce the experimental effort for the non-random two-liquid parameter estimation.
\smallskip

In future work, we plan to demonstrate the benefits of the proposed sequential design methodolgy also on other application examples, especially, on systems with liquid-liquid equilibria and on systems with chemical reactions.  
We also plan to extend the proposed sequential design methodolgy to experimental design problems with constraints on the designs, 
like cost or location constraints. In practice, such constraints are highly relevant, but are rarely taken into account systematically. 
And finally, we will investigate how experimental design 
can further be improved with the help of other design criteria than the linearization-based design criteria used in this work. In particular, it will be interesting to see how and to what extent experimental design for nonlinear models improves when the designs are computed using suitable sampling-based approximations to the parameter uncertainties or to the prediction uncertainties, as opposed to the successive linearization-based approximations used in this work. 

\section*{Acknowledgments}

We gratefully acknowledge funding from the Deutsche Forschungsgemeinschaft (DFG, German Research Foundation) – project number 466397921 – 
within the Priority Programme ``SPP 2331: Machine Learning in Chemical Engineering''.

\section*{Author contributions}

M. Bubel: conceptualization, software, validation, data curation, visualization, writing - review and editing. J. Schmid: conceptualization, methodology, formal analysis, validation, writing - original draft, writing - review and editing. V. Kozachynskyi: conceptualization, software, investigation, validation, resources, data curation, visualization, writing - review and editing. E. Esche: conceptualization, investigation, resources, writing - review and editing, supervision, project administration, funding acquisition. M. Bortz: conceptualization, writing - review and editing, supervision, project administration, funding acquisition.

\begin{small}

\end{small}

\end{document}